\documentclass[12pt]{amsart}
\usepackage[all]{xypic}
\usepackage{amsthm}
\usepackage{amssymb}
\usepackage{verbatim}
\usepackage{enumerate}
\usepackage{graphicx}
\title{More on setwise climbability properties}
\date{}
\author{Bernhard K\"onig and Yasuo~Yoshinobu}
\thanks{The second author is partially supported by JSPS KAKENHI Grant Number 18K03394.}

\address{(Bernhard K\"onig)\newline
Risk Management Function\newline
Allianz Versicherungs-AG\newline
K\"oniginstr. 28, 80802 M\"unchen\newline
Germany}
\email{bkoenig.smtp@gmail.com}

\address{(Yasuo Yoshinobu)\newline
Graduate School of Informatics\newline
Nagoya University\newline
Furo-cho, Chikusa-ku, Nagoya 464-8601\newline
Japan}
\email{yosinobu@i.nagoya-u.ac.jp}

\theoremstyle{definition}
\newtheorem{dfn}{Definition}[section]
\newtheorem{prop}[dfn]{Proposition}
\newtheorem{thm}[dfn]{Theorem}
\newtheorem{lma}[dfn]{Lemma}
\newtheorem{cor}[dfn]{Corollary}

\newtheorem{qtn}[dfn]{Question}

\newcommand{\restrict}{\upharpoonright}

\newcommand{\force}{\Vdash}

\newcommand{\concat}[2]{{{#1}^\smallfrown#2}}
\newcommand{\seq}[1]{{\langle#1\rangle}}

\newcommand{\rmi}{\mathrm{I}}
\newcommand{\rmii}{\mathrm{II}}

\makeatletter
\renewcommand{\p@enumii}{}
\makeatother
\makeatletter
\@namedef{subjclassname@2020}{\textup{2020} Mathematics Subject Classification}
\makeatother
\begin{document}
\begin{abstract}
We introduce two types of variations of {\it setwise climbability properties}, which have been introduced by the second author as fragments of Jensen's square principles. We show that variations of the first type are equivalent to known principles and that they are consistent with the Proper Forcing Axiom ($\mathrm{PFA}$). On the other hand, those of the second type can be characterized as Martin-type axioms for some classes of posets defined in terms of a new variation of generalized Banach-Mazur games, and they are no longer consistent with $\mathrm{PFA}$. We also study how large fragments of $\mathrm{PFA}$ are consistent with these principles.
\end{abstract}
\maketitle

%%%%%%%%%%%%%%%%%%%%%%%%%%%%%%%%%%%%%%%%%%%%
\section{Introduction}\label{sec:introduction}
%%%%%%%%%%%%%%%%%%%%%%%%%%%%%%%%%%%%%%%%%%%%
% Section 0 Introduction

Several preceding studies revealed that Jensen's square principles and some of their fragments can be regarded as Martin-type axioms for posets with some game closure properties. For example, $\square_\kappa$ is equivalent to $\mathrm{MA}_\kappa$ for $(\kappa+1)$-strategically closed posets (\cite{velleman_strategy}, \cite{ishiuyoshinobu02}). Also, those fragments can be separated from each other by observing the extent of their consistency with some well-known forcing axioms, like $\mathrm{PFA}$ or $\mathrm{MA}^+(\text{$\sigma$-closed})$.

The {\it setwise climbability properties\/} $\mathrm{SCL^-}$ and $\mathrm{SCL}$, introduced in \cite{yoshinobu17_starvar}, are two of such fragments of $\square_{\omega_1}$. They can be characterized respectively as $\mathrm{MA}_{\omega_2}$ for a class of posets with some game closure properties defined in terms of a variation (called as the {\it $*$-variation}) of generalized Banach-Mazur games. 

In this paper, we introduce two types of natural variations of the setwise climbability properties, respectively named as the `full' and the `end-extension' variations. The `full' variations, $\mathrm{SCL^-_f}$ and $\mathrm{SCL_f}$, turn out to be equivalent to some (combinations of) existing principles. The `end-extension' variations, $\mathrm{SCL^-_e}$ and $\mathrm{SCL_e}$, are new and characterized in terms of a further variation of the $*$-variation (labelled {\it $**$-variation}) of generalized Banach-Mazur games. Though the difference between the rules of $*$- and $**$-variations seems to be small, corresponding game closure properties are very different: forcing with posets having the game closure property defined from the $*$-variation preserves $\mathrm{PFA}$, whereas forcing with posets having the game closure property defined from the $**$-variation does not. By taking closer look on the $**$-variation, we also found a separation of two seemingly close fragments of $\mathrm{PFA}$.

This paper proceeds as follows. In the rest of this section we define relevant fragments of square principles and game closure properties and review facts about them which are known in the preceding studies. In Section 2 we introduce and analyze the `full' variations of the setwise climbability properties. In Section 3 we introduce and analyze the `end-extension' variations of the principles. We also introduce the $**$-variation of generalized Banach-Mazur games, and characterize the `end-extension' variations using it. In Section 4 we introduce two strengthenings of properness, {\it absolute\/} and {\it indestructible\/} properness, and separate $\mathrm{MA}_{\omega_1}(\text{absolutely proper})$ and $\mathrm{MA}_{\omega_1}(\text{indestructibly proper})$ using the $**$-variation. 
% Section 1 Notations and definitions

Our notations are mostly standard. $S^\beta_\gamma$ denotes the set $\{\alpha<\omega_\beta\mid\operatorname{cf}\alpha=\omega_\gamma\}$. $\mathrm{Lim}$ denotes the class of limit ordinals. Throughout this paper we treat $0$ as a nonlimit ordinal. For a set $S$ of ordinals $\mathrm{o.t.}(S)$ denotes the order type of $S$ and $\mathrm{acc}(S)$ denotes the set of accumulation points of $S$. For a poset $\mathbb{P}$, $\mathcal{B}(\mathbb{P})$ denotes the Boolean completion of $\mathbb{P}$.

\medskip\noindent
\underline{Combinatorial Principles}
\smallskip
\begin{dfn}\label{dfn:properties}
Let $\kappa$ be an uncountable cardinal.
\begin{enumerate}[(1)]
\item\label{item:square}(\cite{finestructure})
A sequence $\seq{C_\alpha\mid\alpha\in\mathrm{Lim}\cap\kappa^+}$ is said to be a {\it $\square_\kappa$-sequence\/} if for every $\alpha\in\mathrm{Lim}\cap\kappa^+$
\begin{enumerate}[(i)]
  \item $C_\alpha$ is a club subset of $\alpha$, 
  \item $\mathrm{o.t.}(C_\alpha)\leq\kappa$, and
  \item $C_\alpha\cap\gamma=C_\gamma$ for every $\gamma\in\mathrm{acc}(C_\alpha)$.
\end{enumerate}
$\square_\kappa$ denotes the statement that there exists a $\square_\kappa$-sequence.
\item\label{item:appr}(\cite{shelah79}, \cite{foreman97})
A sequence $\seq{C_\alpha\mid\alpha\in\mathrm{Lim}\cap\kappa^+}$ is said to be {\it an $\mathrm{AP}_\kappa$-sequence\/} if there exists a club subset $C$ of $\mathrm{Lim}\cap\kappa^+$ such that for every $\alpha\in C$
\begin{enumerate}[(i)]
  \item $C_\alpha$ is a club subset of $\alpha$,
  \item $\mathrm{o.t.}(C_\alpha)\leq\kappa$ and
  \item for every $\gamma<\alpha$ there exists $\beta<\alpha$ such that $C_\alpha\cap\gamma=C_\beta$.
\end{enumerate}
$\mathrm{AP}_\kappa$\footnote{`$\mathrm{AP}$' stands for `approachability'. This property was originally defined in terms of Shelah's approachability ideal, which was introduced in \cite{shelah79}. \cite{foreman97} introduces the `approachability property' in a seemingly stronger form, and states it is equivalent to the present form under $\mathrm{GCH}$. In fact, it is not hard to see they are equivalent without any assumption (see \cite{yoshinobu03:_approac} for proof).} denotes the statement that there exists an $\mathrm{AP}_\kappa$-sequence.
\item\label{item:cp}(\cite{yoshinobu13:_oper}) A function $f:\kappa^+\to\kappa$ is said to be a {\it $\mathrm{CL}_{\kappa}$-function\/}\footnote{`$\mathrm{CP}$' has been used instead of `$\mathrm{CL}$' in the original paper.} if for every $\beta\in\mathrm{Lim}\cap\kappa^+$ there exists $C$ such that
\begin{enumerate}[(i)]
  \item $C$ is a club subset of $\beta$,
  \item $\mathrm{o.t.}(C)\leq\kappa$ and
  \item $f(\alpha)=\mathrm{o.t.}(C\cap\alpha)$ for every $\alpha\in C$.
\end{enumerate}
$\mathrm{CL}_\kappa$ denotes the statement that there exists a $\mathrm{CL}_\kappa$-function.
\end{enumerate}

We say a sequence $\seq{C_\alpha\mid\alpha\in S^2_0}$ is an {\it $S^2_0$-system} if $C_\alpha$ is a countable unbounded subset of $\alpha$ for every $\alpha\in S^2_0$.

\begin{enumerate}[(1)]
\setcounter{enumi}{3}
\item\label{item:scpminus}(\cite{yoshinobu17_starvar}) An $S^2_0$-system $\seq{C_\alpha\mid\alpha\in S^2_0}$ is said to be an {\it $\mathrm{SCL^-}$-system\/}\footnote{`$\mathrm{SCP}$' has been used instead of `$\mathrm{SCL}$' in the original paper.} if for every $\beta\in S^2_1$ there exists $C$ such that
\begin{enumerate}[(i)]
  \item $C$ is a club subset of $\beta\cap S^2_0$,
  \item $\mathrm{o.t.}(C)=\omega_1$ and
  \item $\seq{C_\alpha\mid\alpha\in C}$ is $\subseteq$-increasing and continuous.
\end{enumerate}
$\mathrm{SCL}^-$ denotes the statement that there exists an $\mathrm{SCL}^-$-system.
\item\label{item:scp}(\cite{yoshinobu17_starvar}) A pair $\seq{f, \seq{C_\alpha\mid\alpha\in S^2_0}}$ of a function $f:S^2_0\to\omega_1$ and an $S^2_0$-system $\seq{C_\alpha\mid\alpha\in S^2_0}$ is said to be an {\it $\mathrm{SCL}$-pair\/} if for every $\beta\in S^2_1$ there exists $C$ such that
\begin{enumerate}[(i)]
  \item $C$ is a club subset of $\beta\cap S^2_0$,
  \item $\mathrm{o.t.}(C)=\omega_1$,
  \item $f(\alpha)=\mathrm{o.t.}(C\cap\alpha)$ for every $\alpha\in C$ and
  \item $\seq{C_\alpha\mid\alpha\in C}$ is $\subseteq$-increasing and continuous.
\end{enumerate}
$\mathrm{SCL}$ denotes the statement that there exists an $\mathrm{SCL}$-pair.
\end{enumerate}
\end{dfn}

%We use the following notations.
%\begin{eqnarray*}
%Q&=&\{\omega_1\alpha\mid\alpha\in S^2_0\}\\
%R_0&=&\{\omega_1\alpha+\omega\gamma\mid\alpha\in S^2_0\land 1\leq\gamma<\omega_1\}\\
%R_1&=&\{\omega_1\alpha+\omega\gamma\mid\alpha\in \omega_2\setminus S^2_0\land 1\leq\gamma<\omega_1\}\\
%T&=&\{\omega_1\beta\mid\beta\in\mathrm{Succ}\cap\omega_2\}\\
%U&=&\{\omega_1\beta\mid\beta\in S^2_1\}.
%\end{eqnarray*}

\smallskip\noindent
\underline{Generalized Banach-Mazur games}
\begin{dfn}\normalfont\label{dfn:bmgame}
Let $\mathbb{P}$ be a separative poset. We let $\mathbb{P}^*$ denote the set of at most countable compatible subsets of $\mathbb{P}$.
\begin{enumerate}[(1)]
\item\label{item:foreman}(\cite{foreman83:_games_boolean})
For an ordinal $\lambda$, $G_\lambda(\mathbb{P})$ denotes the following two-player game, played by Players $\rmi$ and $\rmii$. The game develops in (at most) $\lambda$ rounds. In the $\alpha$-th round, for each nonlimit $\alpha<\lambda$, Player $\rmi$ is first required to choose a move $a_\alpha\in\mathbb{P}$, stronger than $b_\beta$ in cases $\alpha=\beta+1$, and then Player $\rmii$ is required to choose $b_\alpha\in\mathbb{P}$ stronger than $a_\alpha$. In the $\alpha$-th round for limit $\alpha<\lambda$, only Player $\mathrm{II}$ is required to choose $b_\alpha$ stronger than all preceding moves. Then the game develops in the following way.
$$
\begin{matrix}
\text{Player}\ \rmi: & a_0\phantom{b_0} &  a_1\phantom{b_1} & \cdots & \phantom{b_\omega} & a_{\omega+1}\phantom{b_{\omega+1}} & \cdots \\
\text{Player}\ \rmii: & \phantom{a_0}b_0 & \phantom{a_1}b_1 & \cdots & b_\omega &  \phantom{a_{\omega+1}}b_{\omega+1} & \cdots
\end{matrix}
$$
Player $\mathrm{II}$ wins the game if she was able to make legal moves in all rounds.

$G^\rmi_\lambda(\mathbb{P})$ denotes a similar game to $G_\lambda(\mathbb{P})$, but in this game both players are required to make stronger and stronger moves in all rounds, that is, the game develops in the following way.
$$
\begin{matrix}
\text{Player}\ \rmi: & a_0\phantom{b_0} &  a_1\phantom{b_1} & \cdots & a_\omega\phantom{b_\omega} & a_{\omega+1}\phantom{b_{\omega+1}} & \cdots \\
\text{Player}\ \rmii: & \phantom{a_0}b_0 & \phantom{a_1}b_1 & \cdots & \phantom{a_\omega}b_\omega &  \phantom{a_{\omega+1}}b_{\omega+1} & \cdots
\end{matrix}
$$
Player $\rmii$ wins the game if she was able to make moves in all rounds, without driving Player $\rmi$ unable to make legal moves on the way.

A funtcion $\sigma:X\to\mathbb{P}$ for some $X\subseteq\bigcup_{\alpha<\lambda}{}^\alpha\mathbb{P}$ is called a {\it strategy\/} for $\mathbb{P}$.

In a play of $G_\lambda(\mathbb{P})$ or $G^\rmi_\lambda(\mathbb{P})$, Player $\rmii$ is said to {\it play by $\sigma$} if in the $\alpha$-th round for each $\alpha<\lambda$ she chooses $\sigma(\vec{a})$ as her move as long as possible, where $\vec{a}$ denotes the sequence consisting of preceding moves of Player $\rmi$.

$\mathbb{P}$ is said to be {\it $\lambda$-strategically closed\/} ({\it resp.} {\it $\lambda$-strongly strategically closed\/}) if there exists a strategy $\sigma$ such that Player $\rmii$ wins any play of $G_\lambda(\mathbb{P})$ ({\it resp.} $G^\rmi_\lambda(\mathbb{P})$) as long as she plays by $\sigma$.

\item\label{item:tacoper}(\cite{yoshinobu13:_oper})
A function $\sigma: X\to\mathbb{P}$ is said to be an {\it operation\/} ({\it resp.} a {\it tactic}) if $X\subseteq\mathrm{On}\times\mathcal{B}(\mathbb{P})$ ({\it resp.} $X\subseteq\mathcal{B}(\mathbb{P})$).

In a play of $G_\lambda(\mathbb{P})$ for an ordinal $\lambda$, Player $\rmii$ is said to {\it play by $\sigma$} if in the $\alpha$-th round for each $\alpha<\lambda$ she chooses $\sigma(\alpha, a)$ ({\it resp.} $\sigma(a)$) as her move as long as possible, where $a$ denotes the boolean infimum of preceding moves of Player $\rmi$.

$\mathbb{P}$ is said to be {\it $\lambda$-operationally closed\/} ({\it resp.} {\it $\lambda$-tactically closed\/}) if there exists an operation ({\it resp.} a tactic) $\sigma$ such that Player $\rmii$ wins any play of $G_\lambda(\mathbb{P})$ as long as she plays by $\sigma$.

\item\label{item:stargame}(\cite{yoshinobu17_starvar}) $G^*(\mathbb{P})$ denotes the following variation of $G_{\omega_1+1}(\mathbb{P})$. In the $\alpha$-th round for nonlimit $\alpha<\omega_1$, Player $\rmi$ is first required to choose a move $A_\alpha\in\mathbb{P}^*$ containing all preceding moves made by himself, so that the boolean infimum $a_\alpha$ of $A_\alpha$ computed in $\mathcal{B}(\mathbb{P})$ is stronger than $b_\beta$ in cases $\alpha=\beta+1$, and then Player $\rmii$ is required to choose a move $b_\alpha\in\mathbb{P}$ stronger than $a_\alpha$ (equivalently, stronger than every condition in $A_\alpha$). In the $\alpha$-th round for limit $\alpha<\omega_1$, Player $\rmi$ is forced to choose $A_\alpha=\bigcup_{\gamma<\alpha}A_\gamma$ as his move, and then Player $\rmii$ is required to choose $b_\alpha$ stronger than every condition in $A_\alpha$. In the $\omega_1$-th round only Player $\rmii$ is required to choose $b_{\omega_1}$ stronger than all of her preceding moves. Then the game develops in the following way.
$$
\begin{matrix}
\text{Player}\ \rmi: & A_0\phantom{b_0} &  A_1\phantom{b_1} & \cdots & A_\omega\phantom{b_\omega} & A_{\omega+1}\phantom{b_{\omega+1}} & \cdots \\
\text{Player}\ \rmii: & \phantom{A_0}b_0 & \phantom{A_1}b_1 & \cdots & \phantom{A_\omega}b_\omega &  \phantom{A_{\omega+1}}b_{\omega+1} & \cdots
\end{matrix}
$$
Player $\mathrm{II}$ wins the game if she was able to make legal moves in all rounds (she is responsible for keeping $A_\alpha\in\mathbb{P}^*$ in limit rounds).

A function $\sigma:X\to\mathbb{P}$ is said to be a {\it $*$-operation} ({\it resp.} a {\it $*$-tactic}) for $\mathbb{P}$ if $X\subseteq\mathrm{On}\times\mathbb{P}^*$ ({\it resp.} $X\subseteq\mathbb{P}^*$).

In a play of $G^*(\mathbb{P})$, Player $\rmii$ is said to {\it play by $\sigma$} if in the $\alpha$-th round for each $\alpha<\omega_1$ she chooses $\sigma(\alpha, A_\alpha)$ ({\it resp.} $\sigma(A_\alpha)$) as her move as long as it is possible, where $A_\alpha$ denotes Player $\rmi$'s move in the $\alpha$-th round.

$\mathbb{P}$ is said to be {\it $*$-operationally closed\/} ({\it resp.} {\it $*$-tactically closed\/} if there exists a $*$-operation ({\it resp.} a $*$-tactic) $\sigma$ such that Player $\rmii$ wins any play of $G^*(\mathbb{P})$ as long as she plays by $\sigma$ (that is, she has a legal move even in the $\omega_1$-th round).
\end{enumerate}
\end{dfn}

\smallskip\noindent
\underline{Martin-type axioms}

\smallskip
Martin-type axioms appear in this article in two different ways: Firstly, several fragments of Jensen's square principles are understood as $\mathrm{MA}_\lambda$ (mostly $\lambda\geq\omega_2$) for some class of posets with some `game-closure' properties. Secondly, we use well-known $\mathrm{MA}_{\omega_1}$-type axioms like $\mathrm{MM}$ or $\mathrm{PFA}$, to compare various class of posets with such game-closure properties, by observing the extent forcing with posets in those classes preserve or destroy these axioms.

\begin{dfn}\normalfont
Let $\mathbb{P}$ be a poset and $\lambda$ a cardinal. $\mathrm{MA}_\lambda(\mathbb{P})$ denotes the following statement: For any $p\in\mathbb{P}$ and any family of dense subsets $\{D_\alpha\mid\alpha<\lambda\}$, there exists a filter $\mathcal{F}\subseteq\mathbb{P}$ such that $p\in\mathcal{F}$ and $D_\alpha\cap\mathcal{F}\not=\emptyset$ for all $\alpha<\lambda$.

$\mathrm{MA}^+(\mathbb{P})$ denotes the following statement: For any $p\in\mathbb{P}$, any $\mathbb{P}$-name $\dot{S}$ for a subset of $\omega_1$ and any family of dense subsets $\{D_\alpha\mid\alpha<\omega_1\}$, there exists a filter $\mathcal{F}\subseteq\mathbb{P}$ such that $p\in\mathcal{F}$, $D_\alpha\cap\mathcal{F}\not=\emptyset$ for all $\alpha<\omega_1$ and that $\dot{S}_\mathcal{F}=\{\alpha\in\omega_1\mid\exists q\in\mathcal{F}(q\force_\mathbb{P}\check{\alpha}\in\dot{S})\}$ is stationary in $\omega_1$.

For a class $\Gamma$ of posets and an ordinal $\lambda$, $\mathrm{MA}_\lambda(\Gamma)$ ({\it resp.} $\mathrm{MA}^+(\Gamma)$) denotes the statement that $\mathrm{MA}_\lambda(\mathbb{P})$ ({\it resp.} $\mathrm{MA}^+(\mathbb{P})$) holds for all $\mathbb{P}\in\Gamma$.

$\mathrm{MM}$ denotes $\mathrm{MA}_{\omega_1}(\text{preserving stationary subsets of $\omega_1$})$ and $\mathrm{PFA}$ denotes $\mathrm{MA}_{\omega_1}(\text{proper})$.
\end{dfn}

\smallskip\noindent
\underline{Known Facts}

%$$
%\xymatrix{
%\square_\kappa \ar@{=>}[r] \ar@{=>}[d] & \mathrm{AP}_\kappa & & \square_{\omega_1} \ar@{=>}[d] \ar@{=>}[rr] & & \mathrm{AP}_{\omega_1}\\
%\mathrm{CP}_\kappa & & & \mathrm{SCP} \ar@{=>}[ld] \ar@{=>}[rd] & & \\
%& & \mathrm{SCP}^- & & \mathrm{CP}_{\omega_1} &
%}
%$$
%
%$$
%\xymatrix{
%%\text{$\omega_2$-directed closed} \ar@{=>}[d]\\
%\text{$\omega_2$-closed} \ar@{=>}[d]\\
%\text{$(\omega_1+1)$-tactically closed}\\
%\text{$(\omega_1+1)$-operationally closed} \text{$*$-tactically closed}\\
%\text{$*$-operationally closed} \text{$(\omega_1+1)$-strongly strategically closed}\\
%\text{$(\omega_1+1)$-strategically closed}
%}
%$$

\smallskip
We list facts about combinatorial principles and game closure properties raised above that are easily observed or known in preceding studies.

Let us first observe easy implications between combinatorial principles and game closure properties respectively.

\begin{prop}\normalfont
Let $\kappa$ be any uncountable cardinal.
\begin{enumerate}[(1)]
\item $\square_\kappa$ implies both $\mathrm{AP}_\kappa$ and $\mathrm{CL}_\kappa$.
\item $\square_{\omega_1}$ implies $\mathrm{SCL}$.
\item $\mathrm{SCL}$ implies both $\mathrm{SCL}^-$ and $\mathrm{CL}_{\omega_1}$.
\end{enumerate}
\end{prop}

\begin{prop}\normalfont
Let $\kappa$ be an infinite cardinal.
\begin{enumerate}[(1)]
\item Every $\kappa^+$-closed poset is both $(\kappa+1)$-tactically closed and $(\kappa+1)$-strongly strategically closed.
\item Every $(\kappa+1)$-tactically closed poset is $(\kappa+1)$-operationally closed.
\item Every $(\kappa+1)$-operationally closed poset and every $(\kappa+1)$-strongly strategically closed is $(\kappa+1)$-strategically closed.
\item Every $(\omega_1+1)$-tactically closed poset is $*$-tactically closed.
\item Every $(\omega_1+1)$-operationally closed poset and every $*$-tactically closed poset is $*$-operationally closed.
\end{enumerate}
\end{prop}

The combinatorial principles and game closure properties raised above are related in the following way.
\begin{prop}\normalfont\label{prop:naturalposet}
Let $\kappa$ be an uncountable cardinal.
\begin{enumerate}[(1)]
\item\label{item:psquare} The poset $\mathbb{P}_{\square_\kappa}$, which naturally adds a $\square_\kappa$-sequence by approximation, is $(\kappa+1)$-strategically closed.
\item The poset $\mathbb{P}_{\mathrm{AP}_\kappa}$, which naturally adds an $\mathrm{AP}_\kappa$-sequence by approximation, is $(\kappa+1)$-strongly strategically closed.
\item The poset $\mathbb{P}_{\mathrm{CL}_\kappa}$, which naturally adds a $\mathrm{CL}_\kappa$-function by approximation, is $(\kappa+1)$-operationally closed.
\item The poset $\mathbb{P}_{\mathrm{SCL^-}}$, which naturally adds an $\mathrm{SCL^-}$-system by approximation, is $*$-tactically closed.
\item\label{item:PSCL} The poset $\mathbb{P}_{\mathrm{SCL}}$, which naturally adds an $\mathrm{SCL}$-pair by approximation, is $*$-operationally closed.
\end{enumerate}
\end{prop}

Using this, one can characterize these principles as Martin-type axioms.

\begin{thm}\normalfont\label{thm:squareasmartin}
Let $\kappa$ be an uncountable cardinal.
\begin{enumerate}[(1)]
\item\label{item:squarekappa}(\cite{ishiuyoshinobu02}) $\square_\kappa$ is equivalent to $\mathrm{MA}_{\kappa^+}(\text{$(\kappa+1)$-strategically closed})$.
\item\label{item:apkappa}(\cite{yoshinobu03:_approac}) $\mathrm{AP}_\kappa$ is equivalent to $\mathrm{MA}_{\kappa^+}(\text{$(\kappa+1)$-strongly strategically closed})$.
\item(\cite{yoshinobu13:_oper}) $\mathrm{CL}_\kappa$ is equivalent to $\mathrm{MA}_{\kappa^+}(\text{$(\kappa+1)$-operationally closed})$.
\item(\cite{yoshinobu17_starvar}) $\mathrm{SCL}^-$ is equivalent to $\mathrm{MA}_{\omega_2}(\text{$*$-tactically closed})$.
\item(\cite{yoshinobu17_starvar}) $\mathrm{SCL}$ is equivalent to $\mathrm{MA}_{\omega_2}(\text{$*$-operationally closed})$.
\end{enumerate}
\end{thm}

\noindent\underline{Remark} As for (\ref{item:squarekappa}) above, more precisely, in \cite{ishiuyoshinobu02} it was proved that $\square_\kappa$ is equivalent to the statement that every $(\kappa+1)$-strategically closed poset is $\kappa^+$-strategically closed, the property for which $\mathrm{MA}_{\kappa^+}$ immediately follows. The other direction is obtained by applying $\mathrm{MA}_{\kappa^+}$ to $\mathbb{P}_{\square_\kappa}$ which is $(\kappa+1)$-strategically closed as stated in Proposition \ref{prop:naturalposet}(\ref{item:psquare}). The role of \cite{yoshinobu03:_approac} in the proof of (\ref{item:apkappa}) is similar.

Finally, these principles (and these game closure properties) can be separated by observing their consistency with (or the extent of preservation of) well-known forcing axioms.
\begin{thm}\normalfont\label{thm:separate}
\begin{enumerate}[(1)]
\item(\cite{koenig04:_fragm_maxim}) Both $\mathrm{PFA}$ and $\mathrm{MA}^+(\text{$\sigma$-closed})$ are preserved under any $\omega_2$-closed forcing.
\item(\cite{koenig04:_fragm_maxim}) $\mathrm{MA}^+(\text{$\sigma$-closed})$ is preserved under any $(\omega_1+1)$-strongly strategically closed forcing.
\item\label{item:pfanegatesap}(independently proved by Foreman and Todorcevic, see \cite{koenig04:_fragm_maxim}) $\mathrm{PFA}$ negates $\mathrm{AP}_{\omega_1}$, and thus may be destroyed by some $(\omega_1+1)$-strongly strategically closed forcing.
\item\label{item:pfapreservedstar}(\cite{yoshinobu17_starvar}) $\mathrm{PFA}$ is preserved under any $*$-operationally closed forcing (and therefore $\mathrm{PFA}$ is consistent with $\mathrm{SCL}$, by Proposition \ref{prop:naturalposet}(\ref{item:PSCL})).
\item(\cite{yoshinobu13:_oper}, \cite{yoshinobu17_starvar}) $\mathrm{MA}^+(\text{$\sigma$-closed})$ negates both $\mathrm{CL}_{\omega_1}$ and $\mathrm{SCL}^-$, and thus may be destroyed both by some $(\omega_1+1)$-operationally closed forcing and by some $*$-tactically closed forcing.
\item(\cite{yoshinobu17_starvar}) Under $\mathrm{MA}^+(\text{$\sigma$-closed})$, no $(\omega_1+1)$-operationally closed forcing forces $\mathrm{SCL}^-$ and no $*$-tactically closed forcing forces $\mathrm{CL}_{\omega_1}$. Thus neither $\mathrm{SCL}^-$ nor $\mathrm{CL}_{\omega_1}$ implies the other.
\end{enumerate}
\end{thm}

\section{The \lq\lq full\rq\rq variations}\label{sec:full}
%%%%%%%%%%%%%%%%%%%%%%%%%%%%%%%%%%%%%%%%%%%%
%%%%%%%%%%%%%%%%%%%%%%%%%%%%%%%%%%%%%%%%%%%%
%\section{The "full" variation}
%%%%%%%%%%%%%%%%%%%%%%%%%%%%%%%%%%%%%%%%%%%%
In this section we analyze the following couple of variations respectively of $\mathrm{SCL^-}$ and $\mathrm{SCL}$.

\begin{dfn}\label{dfn:fullvar}
\begin{enumerate}[(1)]
\item\label{item:scpminusf} An $S^2_0$-system $\seq{C_\alpha\mid\alpha\in S^2_0}$ is said to be an {\it $\mathrm{SCL^-_f}$-system\/} if for every $\beta\in S^2_1$ there exists $C$ such that
\begin{enumerate}[(i)]
  \item $C$ is a club subset of $\beta\cap S^2_0$,
  \item $\mathrm{o.t.}(C)=\omega_1$,
  \item $\seq{C_\alpha\mid\alpha\in C}$ is $\subseteq$-increasing and continuous and
  \item (fullness) $\bigcup_{\alpha\in C}C_\alpha=\beta$.
\end{enumerate}
$\mathrm{SCP^-_f}$ denotes the statement that there exists an $\mathrm{SCL^-_f}$-system.

\item\label{item:scpf} A pair $\seq{f, \seq{C_\alpha\mid\alpha\in S^2_0}}$ of a function $f:S^2_0\to\omega_1$ and an $S^2_0$-system $\seq{C_\alpha\mid\alpha\in S^2_0}$ is said to be an {\it $\mathrm{SCL_f}$-pair\/} if for every $\beta\in S^2_1$ there exists $C$ such that 
\begin{enumerate}[(i)]
  \item $C$ is a club subset of $\beta\cap S^2_0$,
  \item $\mathrm{o.t.}(C)=\omega_1$,
  \item $f(\alpha)=\mathrm{o.t.}(C\cap\alpha)$ for every $\alpha\in C$, 
  \item $\seq{C_\alpha\mid\alpha\in C}$ is $\subseteq$-increasing and continuous and
  \item (fullness) $\bigcup_{\alpha\in C}C_\alpha=\beta$.
\end{enumerate}
$\mathrm{SCL_f}$ denotes the statement that there exists an $\mathrm{SCL_f}$-pair.
\end{enumerate}
\end{dfn}

Note that the only difference between $\mathrm{SCL}^-$ and $\mathrm{SCL^-_f}$ (and that between $\mathrm{SCL}$ and $\mathrm{SCL_f}$) is that the latter requires the fullness condition.

It turns out that these variations are equivalent to some (combination of) known principles.

\begin{thm}\normalfont\label{thm:scpf}
\begin{enumerate}[(1)]
\item\label{item:scpminusfequiv} $\mathrm{SCL^-_f}$ is equivalent to $\mathrm{SCL}^-+\mathrm{CL}_{\omega_1}$.
\item\label{item:scpequiv} $\mathrm{SCL_f}$ is equivalent to $\mathrm{SCL}$.
\end{enumerate}
\end{thm}

To prove this theorem we will use the following lemma.

\begin{dfn}\label{dfn:cpflat}
A function $f:S^2_0\to\omega_1$ is said to be a {\it $\mathrm{CL}^\flat$-function\/} if for every $\beta\in S^2_1$ there exists $C$ such that
\begin{enumerate}[(i)]
  \item $C$ is a club subset of $\beta\cap S^2_0$,
  \item $\mathrm{o.t.}(C)=\omega_1$ and
  \item $f$ is continuous and strictly increasing on $C$.
\end{enumerate}
$\mathrm{CL}^\flat$ denotes the statement that there exists a $\mathrm{CL}^\flat$-function.
\end{dfn}

\begin{lma}\label{lma:flattoregular}
$\mathrm{CL}^\flat$ is equivalent to $\mathrm{CL}_{\omega_1}$.
\end{lma}
\proof It is clear that $\mathrm{CL}_{\omega_1}$ implies $\mathrm{CL}^\flat$. To see the other direction, suppose $f:S^2_0\to\omega_1$ is a $\mathrm{CL}^\flat$-function. We define $g:\omega_2\to\omega_1$ by
$$
\begin{cases}
g(\omega\xi+1+n)&=n\quad\text{for $\xi<\omega_2$ and $n<\omega$,}\\
g(\omega_1\alpha)&=f(\alpha)\quad\text{for $\alpha\in S^2_0$,}\\
g(\omega_1\alpha+\omega(1+\gamma))&=\gamma\quad\text{for $\alpha<\omega_2$ and $\gamma<\omega_1$ and}\\
g(\xi)&=0\quad\text{for $\xi=0$ or $\xi\in S^2_1$}.
\end{cases}
$$
The first line in the above definition deals with the successor ordinals, the second line the multiples of $\omega_1$ in $S^2_0$, the third line the other ordinals in $S^2_0$ and the last line the other ordinals below $\omega_2$. We will show that $g$ is a $\mathrm{CL}_{\omega_1}$-function. Let first $\beta\in S^2_0$. If $\beta$ is of the form $\beta=\omega(\xi+1)$ for some $\xi<\omega_2$, let $$C=\{\omega\xi+1+n\mid n<\omega\}.$$ Otherwise, we may pick an increasing sequence $\seq{\alpha_n|n<\omega}$ of limit ordinals converging to $\beta$. Then let $$C=\{\omega\alpha_n+1+n\mid n<\omega\}.$$ In both cases it is easy to see that $C$ is a club subset of $\beta$ and $g$ maps $C$ order-isomorphically onto $\omega$.

Now let $\beta\in S^2_1$. Note that $\beta$ is of the form $\omega_1\delta$, where $\delta$ is a successor ordinal or $\delta\in S^2_1$. In the case $\delta=\alpha+1$, let
$$
C=\{\omega_1\alpha+\omega(1+\gamma)\mid\gamma<\omega_1\}.
$$
Clearly $C$ is a club subset of $\beta$. Note that $C$ consists only of ordinals in $S^2_0$ which are not multiples of $\omega_1$, and thus it is easy to see that $g$ maps $C$ order-isomorphically onto $\omega_1$.

Now suppose $\delta\in S^2_1$. Since $f$ is a $\mathrm{CL}^\flat$-function, we can pick a club subset $D$ of $\delta\cap S^2_0$ with $\mathrm{o.t.}(D)=\omega_1$ so that $f$ is continuous and strictly increasing on $D$. Now let $\{\alpha_\eta\mid\eta<\omega_1\}$ be the increasing enumeration of $D$ and let
\begin{eqnarray*}
C&=&\{\omega_1\alpha_0+\omega(1+\gamma)\mid\gamma\leq f(\alpha_0)\}\\
&&\cup\bigcup_{\eta<\omega_1}\{\omega_1\alpha_{\eta+1}+\omega(1+\gamma)\mid f(\alpha_\eta)<\gamma\leq f(\alpha_{\eta+1})\}\\
&&\cup\{\omega_1\alpha_\eta\mid\eta\in \omega_1\cap\mathrm{Lim}\}.
\end{eqnarray*}
Then $C$ is a club subset of $S^2_0$ of order type $\omega_1$. Moreover, it holds that
$$
\begin{cases}
g(\omega_1\alpha_0+\omega(1+\gamma))=\gamma&\text{for $\gamma\leq f(\alpha_0)$},\\
g(\omega_1\alpha_{\eta+1}+\omega(1+\gamma)=\gamma&\text{for $f(\alpha_\eta)<\gamma\leq f(\alpha_{\eta+1})$ and}\\
g(\omega_1\alpha_\eta)=f(\alpha_\eta)&\text{for $\eta\in\omega_1\cap\mathrm{Lim}$.}
\end{cases}
$$
This shows that $g$ maps $C$ order-isomorphically onto $\omega_1$. This completes the proof of the claim that $g$ is a $\mathrm{CL}_{\omega_1}$-function.\qed

%
%nd pick a club subset of $C$ of $\beta\cap S^2_0$ of order type $\omega_1$ so that $f$ is continuous and strictly increasing on $C$. Let $\seq{\alpha_{\xi_\gamma}}_{\gamma<\omega_1}$ be the increasing enumeration of $C$. Let
%$$
%  D:=\bigcup_{\gamma<\omega_1}I_\gamma, \ \text{where}
%$$
%\begin{eqnarray*}
%I_0&:=&[\zeta_{\xi_0}+1, \zeta_{\xi_0}+1+f(\alpha_{\xi_0})],\\
%I_{\gamma+1}&:=&[\zeta_{\xi_{\gamma+1}}+1+f(\alpha_{\xi_\gamma})+1, \zeta_{\xi_{\gamma+1}}+1+f(\alpha_{\xi_{\gamma+1}})]\quad(\text{$\gamma<\omega_1$}),\\
%I_{\lambda}&:=&\{\zeta_{\xi_\lambda}\}\ (\text{$\lambda$: limit $<\omega_1$}).
%\end{eqnarray*}
%Then $D$ is a club subset of $\zeta_\beta=\beta$ of order type $\omega_1$. $g$ takes values less than or equal to $f(\alpha_{\xi_0})$ successively on $I_0$, takes values in $[f(\alpha_{\xi_\gamma})+1, f(\alpha_{\xi_{\gamma+1}})]$ successively on $I_{\gamma+1}$ for $\gamma<\omega_1$, and takes $f(\alpha_{\xi_\lambda})$ on $I_\lambda$ for limit $\lambda<\omega_1$.
%This shows that $g$ increasingly enumerates $\omega_1$ on $D$.
%
%Now define $h:S^2_0\to\omega_1$ by
%$$
%h(\alpha_\gamma):=g(\gamma)\quad\text{(for $\gamma<\omega_2$).}
%$$
%Note that, for each $\beta\in S^2_1$ and a club subset $D$ of $\beta$, $\{\alpha_\gamma\mid\gamma\in D\}$ is also a club subset of $\beta$.
%This assures that $h$ is a $\mathrm{CP}^\flat$-function.\qed
%
\smallskip\noindent
(Proof of Theorem \ref{thm:scpf}) For (\ref{item:scpminusfequiv}), first we show that $\mathrm{SCL^-_f}$ implies $\mathrm{SCL}^-$ and $\mathrm{CL}_{\omega_1}$. It is clear that $\mathrm{SCL^-_f}$ implies $\mathrm{SCL}^-$. Suppose $\seq{C_\alpha\mid\alpha\in S^2_0}$ is an $\mathrm{SCL^-_f}$-system. Then define $g:S^2_0\to\omega_1$ by
$$
g(\alpha)=\sup(C_\alpha\cap\omega_1).
$$
For an arbitrary $\beta\in S^2_1$, let $C\subseteq\beta\cap S^2_0$ be as in Definition \ref{dfn:fullvar} (\ref{item:scpminusf}). Then $\seq{C_\alpha\cap\omega_1\mid\alpha\in C}$ forms a continuous $\subseteq$-increasing sequence of countable subsets of $\omega_1$ and $\bigcup_{\alpha\in C}(C_\alpha\cap\omega_1)=\omega_1$. Therefore $g$ is continuous increasing and converging to $\omega_1$ on $C$. So we can pick a club subset $D$ of $C$ so that $g$ is strictly increasing on $D$. This shows that $g$ is a $\mathrm{CL}^\flat$-function, and by Lemma \ref{lma:flattoregular} we have $\mathrm{CL}_{\omega_1}$.

Next we show that $\mathrm{SCL}^-+\mathrm{CL}_{\omega_1}$ implies $\mathrm{SCL^-_f}$. Suppose $\seq{C_\alpha\mid\alpha\in S^2_0}$ is an $\mathrm{SCL}^-$-system and $f:S^2_0\to\omega_1$ is a $\mathrm{CL}_{\omega_1}$-function. For each nonzero $\gamma<\omega_2$ fix a surjection $g_\gamma:\omega_1\to\gamma$. For each $\alpha\in S^2_0$ let
$$
D_\alpha=C_\alpha\cup\{{g_\gamma}(\xi)\mid\gamma\in C_\alpha\land\xi<f(\alpha)\}.
$$
Clearly $D_\alpha$ is a countable unbounded subset of $\alpha$. We will show that $\seq{D_\alpha\mid\alpha\in S^2_0}$ is  an $\mathrm{SCL^-_f}$-system. For each $\beta\in S^2_1$, we can pick a club subset $C$ of $\beta\cap S^2_0$ so that both $\seq{C_\alpha\mid\alpha\in C}$ and $f\restrict C$ are $\subseteq$-continuous increasing. Then it is easy to see that $\seq{D_\alpha\mid\alpha\in C}$ is also $\subseteq$-continuous increasing. Moreover we have %$\displaystyle\bigcup_{\alpha\in C}D_\alpha=\beta$, because
$$
\beta\supseteq\bigcup_{\alpha\in C}D_\alpha\supseteq\{{g_\gamma}(\xi)\mid\gamma\in\bigcup_{\alpha\in C}C_\alpha\land\xi<\omega_1\}=\bigcup\bigcup_{\alpha\in C}C_\alpha=\beta.
$$

For (\ref{item:scpequiv}), it is enough to show that $\mathrm{SCL}$ implies $\mathrm{SCL_f}$, but this can be done exactly in the same way to the proof of $\mathrm{SCL^-_f}$ from $\mathrm{SCL}^-+\mathrm{CL}_{\omega_1}$.\qed

By Theorem \ref{thm:separate}(\ref{item:pfapreservedstar}) we have:
\begin{cor}\normalfont\label{cor:pfascpf}
$\mathrm{SCL_f}$ (and therefore $\mathrm{SCL^-_f}$) is consistent with $\mathrm{PFA}$.\qed
\end{cor}

%We do not know if $\mathrm{SCP}^-_f$ implies $\mathrm{SCP}$.

%%%%%%%%%%%%%%%%%%%%%%%%%%%%%%%%%%%%%%%%%%%%
\section{The \lq\lq end-extension\rq\rq variations}\label{sec:endext}
%%%%%%%%%%%%%%%%%%%%%%%%%%%%%%%%%%%%%%%%%%%%
%%%%%%%%%%%%%%%%%%%%%%%%%%%%%%%%%%%%%%%%%%%%
%\section{The "extension" variation}
%%%%%%%%%%%%%%%%%%%%%%%%%%%%%%%%%%%%%%%%%%%%
In this section we consider several more variations respectively of $\mathrm{SCL}^-$ and $\mathrm{SCL}$.

\begin{dfn}\label{dfn:endvar}
\begin{enumerate}[(1)]
\item\label{item:scpminuse} An $S^2_0$-system $\seq{C_\alpha\mid\alpha\in S^2_0}$ is said to be an {\it $\mathrm{SCL^-_e}$-system\/} if for every $\beta\in S^2_1$ there exists $C$ such that
\begin{enumerate}[(i)]
  \item $C$ is a club subset of $\beta\cap S^2_0$,
  \item $\mathrm{o.t.}(C)=\omega_1$ and
  \item(end-extension) $\seq{C_\alpha\mid\alpha\in C}$ is end-extending.
\end{enumerate}
$\mathrm{SCL^-_e}$ denotes the statement that there exists an $\mathrm{SCL^-_e}$-system.

\item\label{item:scpe} A pair $\seq{f, \seq{C_\alpha\mid\alpha\in S^2_0}}$ of a function $f:S^2_0\to\omega_1$ and an $S^2_0$-system $\seq{C_\alpha\mid\alpha\in S^2_0}$ is said to be an {\it $\mathrm{SCL_e}$-pair\/} if for every $\beta\in S^2_1$ there exists $C$ such that
\begin{enumerate}[(i)]
  \item $C$ is a club subset of $\beta\cap S^2_0$,
  \item $\mathrm{o.t.}(C)=\omega_1$, 
  \item $f(\alpha)=\mathrm{o.t.}(C\cap\alpha)$ for every $\alpha\in C$ and
  \item(end-extension) $\seq{C_\alpha\mid\alpha\in C}$ is end-extending.
\end{enumerate}
$\mathrm{SCL_e}$ denotes the statement that there exists an $\mathrm{SCL_e}$-pair.
\end{enumerate}
\end{dfn}

Note that the only difference between $\mathrm{SCL}^-$ and $\mathrm{SCL^-_e}$ (and that between $\mathrm{SCL}$ and $\mathrm{SCL_e}$) is that in the latter $\seq{C_\alpha\mid\alpha\in C}$ is required to be end-extending, not only to be $\subseteq$-increasing.

As well as $\mathrm{SCL}^-$ and $\mathrm{SCL}$, $\mathrm{SCL^-_e}$ and $\mathrm{SCL_e}$ can be characterized as Martin-type axioms for classes of posets which are defined in terms of another variation of Banach-Mazur games.

\begin{dfn}\label{dfn:doublestargame}
For a separative poset $\mathbb{P}$, $G^{**}(\mathbb{P})$ denotes the following variation of $G^*(\mathbb{P})$. In the $\alpha$-th round for nonlimit $\alpha<\omega_1$, Player $\rmi$ is first required to choose a move $A_\alpha\in\mathbb{P}^*$ containing all preceding moves made by himself, so that in cases $\alpha=\beta+1$, each condition in $A_\alpha$ newly added in the round is stronger than $b_\beta$, and then Player $\rmii$ is required to choose a move $b_\alpha\in\mathbb{P}$ stronger than every condition in $A_\alpha$. In the $\alpha$-th round for limit $\alpha<\omega_1$, Player $\rmi$ is forced to choose $A_\alpha=\bigcup_{\gamma<\alpha}A_\gamma$ and then Player $\rmii$ is required to choose $b_\alpha$ stronger than every condition in $A_\alpha$. In the $\omega_1$-th round only Player $\rmii$ is required to choose $b_{\omega_1}$ stronger than all preceding moves. The game develops as follows:
$$
\begin{matrix}
\text{Player}\ \rmi: & A_0\phantom{b_0} &  A_1\phantom{b_1} & \cdots & A_\omega\phantom{b_\omega} & A_{\omega+1}\phantom{b_{\omega+1}} & \cdots \\
\text{Player}\ \rmii: & \phantom{A_0}b_0 & \phantom{A_1}b_1 & \cdots & \phantom{A_\omega}b_\omega &  \phantom{A_{\omega+1}}b_{\omega+1} & \cdots
\end{matrix}
$$
Player $\rmii$ wins the game if she was able to make legal moves in all rounds, again being responsible for keeping $A_\alpha\in\mathbb{P}^*$ in limit rounds.
\end{dfn}
\begin{dfn}\label{dfn:doublestarclosure}
$\mathbb{P}$ is said to be {\it $**$-operationally closed\/} ({\it resp.} {\it $**$-tactically closed\/}) if there exists a $*$-operation ({\it resp.} a $*$-tactic) $\sigma$ such that Player $\rmii$ wins any play of $G^{**}(\mathbb{P})$ as long as she plays by $\sigma$ (that is, she has a legal move even in the $\omega_1$-th round).
\end{dfn}
It is clear that every $**$-tactically closed poset is $**$-operationally closed, and since $G^{**}(\mathbb{P})$ is a harder game for Player $\rmi$ than $G^*(\mathbb{P})$, every $*$-operationally ({\it resp.} $*$-tactically) closed poset is $**$-operationally ({\it resp.} $**$-tactically) closed.

%\begin{dfn}
%Let $\mathbb{P}$ be a separative poset, and let
%$$
%\mathbb{P}^*=\{A\subseteq\mathbb{P}\mid 1\leq|A|\leq\omega\land\text{$A$ has a common extension in $\mathbb{P}$}\}.
%$$
%\begin{enumerate}[(1)]
%\item A {\it $*$-tactic\/} for $\mathbb{P}$ is a function $\tau:\mathbb{P}^*\to\mathbb{P}$ such that for every $A\in\mathbb{P}^*$ $\tau(A)$ is a common extension of $A$. In a play of $G^{**}(\mathbb{P})$, Player $\rmii$ is said to {\it play with\/} a $*$-tactic $\tau$ if, for each $\alpha<\omega_1$ she chooses $\tau(A_\alpha)$ as her $\alpha$-th move, responding to Player $\rmi$'s $\alpha$-th move $A_\alpha$. A $*$-tactic $\tau$ is said to be a {\it winning\/} one for $G^{**}(\mathbb{P})$ if Player $\rmii$ always wins $G^{**}(\mathbb{P})$ as long as she plays with $\tau$. A poset $\mathbb{P}$ is {\it $**$-tactically closed\/} if there exists a winning $*$-tactic for $G^{**}(\mathbb{P})$.
%\item A {\it $*$-operation\/} for $\mathbb{P}$ is a function $o:\omega_1\times\mathbb{P}^*\to\mathbb{P}$ such that for every $\eta<\omega_1$ and $A\in\mathbb{P}^*$ $o(\eta, A)$ is a common extension of $A$. In a play of $G^{**}(\mathbb{P})$, Player $\rmii$ is said to {\it play with\/} a $*$-operation $o$ if, for each $\alpha<\omega_1$ she chooses $o(\alpha, A_\alpha)$ as her $\alpha$-th move, responding to Player $\rmi$'s $\alpha$-th move $A_\alpha$. The notion of winning $*$-operations for $G^{**}(\mathbb{P})$ and that of $**$-operationally closed posets are defined similarly as above.
%\end{enumerate}
%\end{dfn}

We define natural posets respectively adding an $\mathrm{SCL^-_e}$-system and an $\mathrm{SCL_e}$-pair.

\begin{dfn}
\begin{enumerate}[(1)]
\item A poset $\mathbb{P}_{\mathrm{SCL^-_e}}$ is defined as follows.

\noindent
$p\in\mathbb{P}_{\mathrm{SCL}^-_e}$ iff $p$ is of the form $\seq{C^p_\alpha\mid\alpha\in S^2_0\cap(\gamma^p+1)}$ for some $\gamma^p\in\omega_2\cap\mathrm{Lim}$, and satisfies
\begin{enumerate}[(a)]
  \item for every $\alpha\in S^2_0\cap(\gamma^p+1)$, $C^p_\alpha$ is a countable unbounded subset of $\alpha$, and
  \item for every $\beta\in S^2_1\cap(\gamma^p+1)$, there exists $C$ such that
  \begin{enumerate}[(i)]
    \item $C$ is a club subset of $S^2_0\cap\beta$,
    \item $\mathrm{o.t.}(C)=\omega_1$ and
    \item $\seq{C^p_\alpha\mid\alpha\in C}$ is end-extending.
  \end{enumerate}
\end{enumerate}
For $p$, $q\in\mathbb{P}_{\mathrm{SCL^-_e}}$, $p\leq_{\mathbb{P}_{\mathrm{SCL^-_e}}}q$ iff $q=p\restrict(\gamma^q+1)$.

\item A poset $\mathbb{P}_{\mathrm{SCL_e}}$ is defined as follows.

\noindent
$p\in\mathbb{P}_{\mathrm{SCL_e}}$ iff $p$ is of the form $\seq{f^p, \mathcal{C}^p}$ such that for some $\gamma^p\in\omega_2\cap\mathrm{Lim}$
\begin{enumerate}[(a)]
  \item $f^p:S^2_0\cap(\gamma^p+1)\to\omega_1$,
  \item $\mathcal{C}^p$ is of the form $\seq{C^p_\alpha\mid\alpha\in S^2_0\cap(\gamma^p+1)}$,
  \item for every $\alpha\in S^2_0\cap(\gamma^p+1)$, $C^p_\alpha$ is a countable unbounded subset of $\alpha$ and
  \item for every $\beta\in S^2_1\cap(\gamma^p+1)$, there exists $C$ such that
  \begin{enumerate}[(i)]
    \item $C$ is a club subset of $S^2_0\cap\beta$,
    \item $\mathrm{o.t.}(C)=\omega_1$,
    \item $f^p(\alpha)=\mathrm{o.t.}(C\cap\alpha)$ for every $\alpha\in C$ and
    \item $\seq{C^p_\alpha\mid\alpha\in C}$ is end-extending.
  \end{enumerate}
\end{enumerate}
\end{enumerate}
For $p$, $q\in\mathbb{P}_{\mathrm{SCL_e}}$, $p\leq_{\mathbb{P}_{\mathrm{SCL_e}}}q$ iff $f^q=f^p\restrict(\gamma^q+1)$ and $\mathcal{C}^q=\mathcal{C}^p\restrict(\gamma^q+1)$.
\end{dfn}

\begin{lma}\label{lma:sigmaclosed}
Both $\mathbb{P}_{\mathrm{SCL^-_e}}$ and $\mathbb{P}_{\mathrm{SCL_e}}$ are $\sigma$-closed.
\end{lma}

\proof We prove the theorem only for $\mathbb{P}_{\mathrm{SCL_e}}$. The case of $\mathbb{P}_{\mathrm{SCL^-_e}}$ is easier. Suppose $\seq{p_n\mid n<\omega}$ is a strictly decreasing sequence in $\mathbb{P}_{\mathrm{SCL_e}}$. Then let $\gamma^p=\sup\{\gamma^{p_n}\mid n<\omega\}$ and define $p$ by
$$
\begin{cases}
f^p\restrict\gamma^p=\bigcup\{f^{p_n}\mid n<\omega\}\\
f^p(\gamma^p)=0\\
\mathcal{C}^p\restrict\gamma^p=\bigcup\{\mathcal{C}^{p_n}\mid n<\omega\}\\
C^p_{\gamma^p}=\{\gamma^{p_n}\mid n<\omega\}.
\end{cases}
$$
Then it is easy to check that $p\in\mathbb{P}_{\mathrm{SCL_e}}$ and extends all $p_n$'s.\qed

\begin{lma}\label{lma:poperandtact}
\begin{enumerate}[(1)]
\item\label{item:ptact} $\mathbb{P}_{\mathrm{SCL^-_e}}$ is $**$-tactically closed.
\item\label{item:poper} $\mathbb{P}_{\mathrm{SCL_e}}$ is $**$-operationally closed.
\end{enumerate}
\end{lma}

\proof We will only prove (\ref{item:poper}). (\ref{item:ptact}) is easier. 

For each $\gamma\in S^2_0$, pick a countable unbounded subset $U_\gamma$ of $\gamma$.

For $\eta<\omega_1$ and $A\in(\mathbb{P}_{\mathrm{SCL_e}})^*$, we define $o(\eta, A)$ as follows: If either $A$ has the strongest condition or $\eta$ is nonlimit, let $o(\eta, A)$ be any proper extension of all conditions in $A$. This is possible since $\mathbb{P}_{\mathrm{SCL_e}}$ is $\sigma$-closed. Otherwise, let $\gamma^{o(\eta, A)}=\sup\{\gamma^p\mid p\in A\}$ and define $o(\eta, A)$ by
$$
\begin{cases}
f^{o(\eta, A)}\restrict\gamma^{o(\eta, A)}=\bigcup\{f^p\mid p\in A\}\\
f^{o(\eta, A)}(\gamma^{o(\eta, A)})=\xi\quad\text{(where $\xi$ is such that $\eta=\omega(1+\xi)$)}\\
\mathcal{C}^{o(\eta, A)}\restrict\gamma^{o(\eta, A)}=\bigcup\{\mathcal{C}^p\mid p\in A\}\\
C^{o(\eta, A)}_{\gamma^{o(\eta, A)}}=\{\gamma^p\mid p\in A\}.
\end{cases}
$$
It is easy to check that in any case $o(\eta, A)$ is a $\mathbb{P}_{\mathrm{SCL_e}}$-condition properly extending all conditions in $A$. Note that since $\mathbb{P}_{\mathrm{SCL_e}}$ is $\sigma$-closed, any play of $G^{**}(\mathbb{P}_{\mathrm{SCL_e}})$ never ends before both players make $\omega_1$ moves. Therefore to prove that $o$ is a winning $*$-operation for $G^{**}(\mathbb{P}_{\mathrm{SCL_e}})$, it is enough to show that whenever $\seq{A_\alpha\mid\alpha<\omega_1}$ is a record of moves of Player $\rmi$ in a play of $G^{**}(\mathbb{P}_{\mathrm{SCL_e}})$ in which Player $\rmii$ plays with $o$, $A=\bigcup_{\alpha<\omega_1}A_\alpha$ has a common extension. But now let $\gamma^q=\sup\{\gamma^p\mid p\in A\}$ and define $q$ by
$$
f^q=\bigcup\{f^p\mid p\in A\}\ \text{and}\ \mathcal{C}^q=\bigcup\{\mathcal{C}^p\mid p\in A\}.
$$
To see that $q$ is a $\mathbb{P}_{\mathrm{SCL_e}}$-condition and therefore is a common extension of $A$, note that
$$
C=\{\gamma^{o(\eta, A_\eta)}\mid\eta\in\omega_1\cap\mathrm{Lim}\}
$$
is a club subset of $\gamma^q$ with $\mathrm{o.t.}(C)=\omega_1$ and for each $\eta=\omega(1+\xi)$
$$
f(\gamma^{o(\eta, A_\eta)})=\xi=\mathrm{o.t.}(C\cap\gamma^{o(\eta, A_\eta)})
$$
holds. Moreover we have
$$
C^q_{\gamma^{o(\eta, A_\eta)}}=\{\gamma^p\mid p\in A_\eta\}
$$
and thus $\seq{C^q_\gamma\mid\gamma\in C}$ is end-extending and continuous.\qed

\begin{thm}\label{thm:doublestar}
\begin{enumerate}[(1)]
  \item\label{item:scpeminustact} $\mathrm{SCL^-_e}$ is equivalent to $\mathrm{MA}_{\omega_2}(\text{$**$-tactically closed})$.
  \item\label{item:scpeoper} $\mathrm{SCL_e}$ is equivalent to $\mathrm{MA}_{\omega_2}(\text{$**$-operationally closed})$.
\end{enumerate}
\end{thm}

\proof
%The proof is similar to the similar characterizations of $\mathrm{SCP}^-$ and $\mathrm{SCP}$ (see \cite{}).

We will only show (\ref{item:scpeoper}). (\ref{item:scpeminustact}) is easier.

\smallskip\noindent
(Claim1) For every $\gamma<\omega_2$,
$$
D_\gamma=\{p\in\mathbb{P}_{\mathrm{SCL_e}}\mid\gamma^p\geq\gamma\}
$$
is dense in $\mathbb{P}_{\mathrm{SCL_e}}$.

\proof By straightforward induction on $\gamma<\omega_2$ using the $*$-operation defined in Lemma \ref{lma:poperandtact}.\qed(Claim1)

Now assume $\mathrm{MA}_{\omega_2}(\text{$**$-operationally closed})$. Then there exists a filter $\mathcal{F}$ on $\mathbb{P}_{\mathrm{SCL_e}}$ such that $\mathcal{F}\cap D_\gamma\not=\emptyset$ for every $\gamma<\omega_2$. Then it is easy to see that $q$ defined by
$$
f^q=\bigcup\{f^p\mid p\in\mathcal{F}\}\ \text{and}\ \mathcal{C}^q=\bigcup\{\mathcal{C}^p\mid p\in\mathcal{F}\}
$$

is an $\mathrm{SCL_e}$-pair and therefore we have $\mathrm{SCL_e}$.

Let us prove the other direction. Assume $\mathrm{SCL_e}$ and let $\seq{f, \seq{C_\alpha\mid\alpha\in S^2_0}}$ be an $\mathrm{SCL_e}$-pair. Let $\mathbb{P}$ be any ${**}$-operationally closed poset. Fix a winning $*$-operation $o$ for $G^{**}(\mathbb{P})$. We may assume that $o$ is defined entirely on $\omega_1\times\mathbb{P}^*$, and that $o(\alpha, A)$ is a common extension of the conditions in $A$ for every $\seq{\alpha, A}\in\omega_1\times\mathbb{P}^*$.
To see $\mathrm{MA}_{\omega_2}(\mathbb{P})$, it is enough to show that $\mathbb{P}$ is $\omega_2$-strategically closed. To this end we will use the following lemma:

\begin{lma}\label{lma:ishiu}(Ishiu and Yoshinobu \cite{ishiuyoshinobu02})
Let $\lambda$ be any infinite cardinal. Then there exists a tree ordering $\prec_\lambda$ on the set $N_\lambda$ of nonlimit ordinals below $\lambda$ satisfying
\begin{enumerate}[(a)]
  \item $\seq{N_\lambda, \prec_\lambda}$ is of height $\omega$,
  \item for every $\alpha, \beta\in N_\lambda$, $\alpha\prec_\lambda\beta$ implies $\alpha<\beta$,
  \item for every $\gamma\leq\lambda$ with $\operatorname{cf}\gamma=\omega$, there exists a branch $b$ of $\seq{N_\lambda, \prec_\lambda}$ such that $\sup b=\gamma$.
\end{enumerate}
\end{lma}

Now we describe the winning strategy $\sigma$ for $G_{\omega_2}(\mathbb{P})$. Consider a play of $G_{\omega_2}(\mathbb{P})$ and for each nonlimit $\gamma<\omega_2$ let $a_\gamma$ denote the $\gamma$-th move of Player $\mathrm{I}$ in the play under consideration. At the $\gamma$-th turn of Player $\mathrm{II}$ for $\gamma<\omega_2$, having $\seq{a_\xi\mid\xi\in(\gamma+1)\cap N_{\omega_2}}$ as the preceding moves of the opponent, $\sigma$ suggests Player $\mathrm{II}$ to make her move $b_\gamma$ as follows:
$$
b_\gamma=
\begin{cases}
o(\mathrm{ht}(\gamma, \prec_{\omega_2}), \{a_\xi\mid\xi\preceq_{\omega_2}\gamma\})&\text{if $\gamma\in N_{\omega_2}$,}\\
o(f(\xi), \{a_{\xi+1}\mid\xi\in C_\gamma\})&\text{if $\gamma\in S^2_0$,}\\
\text{any common extension of $\{a_\xi\mid\xi\in\gamma\cap N_{\omega_2}\}$} &\text{if $\gamma\in S^2_1$.}
\end{cases}
$$ 
Let us show that $\sigma$ is a winning strategy for $G_{\omega_2}(\mathbb{P})$. Note first that it is immediate from our assumption on $o$ that $\sigma$ provides a legal move in every round in any play of $G_{\omega_2}(\mathbb{P})$ as long as the play continues.
Therefore to see that $\sigma$ is a winning strategy, it is enough to show that, in any play of $G_{\omega_2}(\mathbb{P})$, whenever Player $\mathrm{II}$ makes first $\gamma$ moves for a limit ordinal $\gamma<\omega_2$ playing with $\sigma$, those moves have a common extension.

\noindent\medskip
\underline{Case 1} $\gamma\in S^2_0$.

Pick a branch $b$ of $\seq{N_{\omega_2}, \prec_{\omega_2}}$ with $\sup b=\gamma$. Let $\seq{\xi_n\mid n<\omega}$ be the enumeration of nodes in $b$ in the order of $\prec_{\omega_2}$. Then note that for each $n<\omega$, Player $\mathrm{II}$'s $\xi_n$-th move $b_{\xi_n}$ is given by
$$
b_{\xi_n}=\sigma(\{a_\eta\mid\eta\in(\xi_n+1)\cap N_{\omega_2}\})=o(n, \{a_{\xi_i}\mid i\leq n\}).
$$
Therefore
$$
\begin{matrix}
\text{Player}\ \rmi: & \{a_{\xi_0}\}\phantom{b_0} &  \{a_{\xi_0,}, a_{\xi_1}\}\phantom{b_1} & \{a_{\xi_0,}, a_{\xi_1}, a_{\xi_2}\}\phantom{b_2} & \cdots\\
\text{Player}\ \rmii: & \phantom{\{a_{\xi_0}\}}b_{\xi_0} & \phantom{\{a_{\xi_0,}, a_{\xi_1}\}}b_{\xi_1} & \phantom{\{a_{\xi_0,}, a_{\xi_1}, a_{\xi_2}\}}b_{\xi_2} & \cdots
\end{matrix}
$$
forms a play of $G^{**}(\mathbb{P})$ where Player $\mathrm{II}$ plays with $o$, and thus the set $\{b_{\xi_n}\mid n<\omega\}$, which is cofinal in $\{b_\eta\mid\eta<\gamma\}$, has a common extension. 

\noindent\medskip
\underline{Case 2} $\gamma\in S^2_1$.

Let $C$ be a club subset of $\gamma\cap S^2_0$ witnessing $\mathrm{SCL_e}$. Let $\seq{\xi_\nu\mid \nu<\omega_1}$ be the increasing enumeration of $C$, Then for each $\nu<\omega_1$, $b_{\xi_\nu}$ is given by
$$
b_{\xi_\nu}=o(\nu, \{a_\eta\mid\eta\in C_{\xi_\nu}\}).
$$
Since $\seq{C_{\xi_\nu}\mid\nu<\omega_1}$ is end-extending and continuous,
$$
\begin{matrix}
\text{Player}\ \rmi: & \{a_\eta\mid\eta\in C_{\xi_0}\}\phantom{b_{\xi_0}} &   \{a_\eta\mid\eta\in C_{\xi_1}\}\phantom{b_{\xi_1}} & \cdots & \{a_\eta\mid\eta\in C_{\xi_\omega}\}\phantom{b_{\xi_\omega}} & \cdots\\
\text{Player}\ \rmii: & \phantom{\{a_\eta\mid\eta\in C_{\xi_0}\}}b_{\xi_0} & \phantom{\{a_\eta\mid\eta\in C_{\xi_1}\}}b_{\xi_1} & \cdots & \phantom{\{a_\eta\mid\eta\in C_{\xi_\omega}\}}b_{\xi_\omega} & \cdots
\end{matrix}
$$
forms a play of $G^{**}(\mathbb{P})$ where Player $\mathrm{II}$ plays with $o$, and thus the set $\{b_{\xi_\nu}\mid \nu<\omega_1\}$, which is cofinal in $\{b_\eta\mid\eta<\gamma\}$, has a common extension. This completes the proof that $\sigma$ is a winning strategy for $G_{\omega_2}(\mathbb{P})$. \qed(Theorem \ref{thm:doublestar})

Now let us observe some implications from $\mathrm{SCL^-_e}$.

\begin{thm}\normalfont\label{thm:scpap}
$\mathrm{SCL^-_e}$ implies $\mathrm{AP}_{\omega_1}$.
\end{thm}

As a corollary, by Theorem \ref{thm:separate}(\ref{item:pfanegatesap}), $**$-tactically closed forcing may destroy $\mathrm{PFA}$, unlike $*$-operationally forcing.

To prove Theorem \ref{thm:scpap}, we use the following lemma.

\begin{lma}\normalfont\label{lma:aptree}(\cite{yoshinobu03:_approac})
For any uncountable cardinal $\kappa$, $\mathrm{AP}_{\kappa}$ is equivalent to the statement that there exists a tree $\seq{\kappa^+, \prec}$ of height $\kappa$ satisfying
\begin{enumerate}[(i)]
\item\label{item:cond1} For every $\alpha$, $\beta\in\kappa^+$, $\alpha\prec\beta$ implies $\alpha<\beta$.
\item\label{item:cond2} For every limit ordinal $\gamma<\kappa^+$, there exists $b\subseteq\kappa^+$ linear and downward closed with respect to $\prec$ such that $\sup b=\gamma$.
\end{enumerate}
\end{lma}

\noindent
(Proof of Theorem \ref{thm:scpap}) Suppose $\seq{C_\alpha\mid\alpha\in S^2_0}$ is an $\mathrm{SCL^-_e}$-system. We will construct a tree $\seq{\omega_2, \prec}$ of height $\omega_1$ satisfying (\ref{item:cond1}) and (\ref{item:cond2}) of Lemma \ref{lma:aptree}.

We define $\prec$ by:
$$
\alpha\prec\beta\Leftrightarrow\alpha\prec_{\omega_2}\beta\ \text{or}\ C_\beta\cap\alpha= C_\alpha,
$$
where $\prec_{\omega_2}$ is as in Lemma \ref{lma:ishiu}. $\prec$ is the disjoint union of $\prec\restrict N_{\omega_2}=\prec_{\omega_2}$ and $\prec\restrict S^2_0$. $\prec\restrict N_{\omega_2}=\prec_{\omega_2}$ is a tree ordering of height $\omega$ by definition. It is clear that $\prec\restrict S^2_0$ is also a tree ordering, and since each $C_\alpha$ is countable, the height of $\prec\restrict S^2_0$ is $\omega_1$. This assures that $\prec$ is a tree ordering on $\omega_2$ of height $\omega_1$. (\ref{item:cond1}) is clear by definition. For $\gamma\in S^2_0$, there is a branch $b\subseteq N_{\omega_2}$ of $\prec_{\omega_2}$ such that  $\sup b=\gamma$ by definition. For $\gamma\in S^2_1$, by the definition of an $\mathrm{SCL^-_e}$-system, there is a club subset $C$ of $\gamma$ such that $C$ is linearly ordered by $\prec$. Therefore $C$ generates a linearly ordered downward closed subset $b$ of $S^2_0$ with respect to $\prec\restrict S^2_0$, and since $\mathrm{o.t.}(C)=\omega_1$, $\sup b=\sup C=\gamma$. This shows (\ref{item:cond2}).\qed(Theorem \ref{thm:scpap})

\smallskip
We end this section by showing somewhat surprising theorem that $\mathrm{SCL^-_e}$ and $\mathrm{SCL_e}$ are in fact equivalent.

\begin{thm}\normalfont\label{thm:minuse}
$\mathrm{SCL^-_e}$ implies $\mathrm{SCL_e}$ (and thus implies $\mathrm{SCL}$).
\end{thm}

\proof Suppose $\seq{C_\alpha\mid\alpha\in S^2_0}$ is an $\mathrm{SCL^-_e}$-system. Put $\delta_\alpha=\mathrm{o.t.}(C_\alpha)$ for $\alpha\in S^2_0$.

Let
\begin{eqnarray*}
Q&=&\{\omega_1\overline{\alpha}+\omega\gamma\mid\overline{\alpha}\in S^2_0, \gamma<\omega_1\},\\
R&=&\{\omega_1\overline{\alpha}+\omega\gamma\mid\overline{\alpha}\in\omega_2\setminus S^2_0\land 1\leq\gamma<\omega_1\},\\
T&=&\{\omega_1(\beta+1)\mid\beta\in\omega_2\}\ \text{and}\\
U&=&\{\omega_1\beta\mid\beta\in S^2_1\}.
\end{eqnarray*}

Note that $S^2_0$ is the disjoint union of $Q$ and $R$, whereas $S^2_1$ is the disjoint union of $T$ and $U$. Define $f:S^2_0\to\omega_1$ by,  for each $\alpha=\omega_1\overline{\alpha}+\omega\gamma\in S^2_0$ ($\gamma<\omega_1$) letting
$$
f(\alpha)=
\begin{cases}
\delta_{\overline{\alpha}}&(\text{if $\overline{\alpha}\in S^2_0$ and $\gamma=0$})\\
\gamma^-&\text{(if $\gamma>0$)},
\end{cases}
$$
where $\gamma^-$ is such that $1+\gamma^-=\gamma$. Define $\seq{\overline{C}_\alpha\mid\alpha\in S^2_0}$ as follows. For $\alpha=\omega_1\overline{\alpha}+\omega\gamma\in S^2_0$ ($\gamma<\omega_1$), let
$$
\overline{C}_\alpha=
\begin{cases}
\{\omega_1\xi\mid\xi\in C_{\overline{\alpha}}\}\cup\displaystyle\bigcup_{\xi\in\mathrm{acc}(C_{\overline{\alpha}})\cap\overline{\alpha}}[\omega_1\xi, \omega_1\xi+\omega\delta_\xi)\cup[\omega_1\overline{\alpha}, \alpha)&\text{($\alpha\in Q$)}\\
[\omega_1\overline{\alpha}, \alpha)&\text{($\alpha\in R$).}
\end{cases}
$$
It is easy to check that $\overline{C}_\alpha$ is countable unbounded subset of $\alpha$.
We will show that $\seq{f, \seq{\overline{C}_\alpha\mid\alpha\in S^2_0}}$ is an $\mathrm{SCL_e}$-pair.

For $\omega_1(\beta+1)\in T$, let $C=\{\omega_1\beta+\omega\gamma\mid 1\leq\gamma<\omega_1\}$. Clearly $C$ is a club subset of $\omega_1(\beta+1)$ and $\mathrm{o.t.}(C)=\omega_1$. For $1\leq\gamma<\omega_1$ we have $f(\omega_1\beta+\omega\gamma)=\gamma^-$. This shows that $f$ maps $C$ order-isomorphically onto $\omega_1$. Moreover we have
$$
\overline{C}_{\omega_1\beta+\omega\gamma}=
\begin{cases}
\overline{C}_{\omega_1\beta}\cup[\omega_1\beta, \omega_1\beta+\omega\gamma)&(\text{if $\beta\in S^2_0$})\\
[\omega_1\beta, \omega_1\beta+\omega\gamma)&(\text{otherwise}).
\end{cases}
$$
This shows that $\seq{\overline{C}_\alpha\mid\alpha\in C}$ is end-extending.

Now for $\omega_1\beta\in U$, pick a club subset $D$ of $\beta\cap S^2_0$ with $\mathrm{o.t.}(D)=\omega_1$ such that $\seq{C_\alpha\mid\alpha\in D}$ is end-extending. Let $\{\alpha_\eta\mid\eta\in\omega_1\}$ be the increasing enumeration of $D$.

Let
\begin{eqnarray*}
C&=&\{\omega_1\alpha_0+\omega\gamma\mid 1\leq\gamma\leq\delta_{\alpha_0}\}\\
&&\cup\bigcup_{\eta\in\omega_1}\{\omega_1\alpha_{\eta+1}+\omega\gamma\mid\delta_{\alpha_\eta}+1\leq\gamma\leq\delta_{\alpha_{\eta+1}}\}\\
&&\cup\{\omega_1\alpha_\eta\mid\eta\in\mathrm{Lim}\cap\omega_1\}.
\end{eqnarray*}
Then $C$ is a club subset of $\omega_1\beta\cap S^2_0$ with $\mathrm{o.t.}(C)=\omega_1$. Now we have
$$
\begin{cases}
f(\omega_1\alpha_0+\omega\gamma)=\gamma^-&\text{for $1\leq\gamma\leq\delta_{\alpha_0}$}\\
f(\omega_1\alpha_{\eta+1}+\omega\gamma)=\gamma&\text{for $\delta_{\alpha_{\eta}}+1\leq\gamma\leq\delta_{\alpha_{\eta+1}}$ and}\\
f(\omega_1\alpha_\eta)=\delta_{\alpha_\eta}&\text{for $\eta\in\mathrm{Lim}\cap\omega_1$}.
\end{cases}
$$
This shows that $f$ maps $C$ order-isomorphically onto $\omega_1$. Moreover, for $\alpha\in C$, 
$$
\overline{C}_\alpha=
\begin{cases}
\displaystyle\{\omega_1\xi\mid\xi\in C_{\alpha_0}\}\cup\bigcup_{\xi\in\mathrm{acc}(C_{\alpha_0})\cap\alpha_0}[\omega_1\xi, \omega_1\xi+\omega\delta_\xi)\cup[\omega_1\alpha_0, \alpha)\\
\quad\text{(for $\alpha=\omega_1\alpha_0+\omega\gamma$, $1\leq\gamma\leq\delta_{\alpha_0}$)}\\
\displaystyle\{\omega_1\xi\mid\xi\in C_{\alpha_{\eta+1}}\}\cup\bigcup_{\xi\in\mathrm{acc}(C_{\alpha_{\eta+1}})\cap\alpha_{\eta+1}}[\omega_1\xi, \omega_1\xi+\omega\delta_\xi)\cup[\omega_1\alpha_{\eta+1}, \alpha)\\
\quad\text{(for $\alpha=\omega_1\alpha_{\eta+1}+\omega\gamma$, $\delta_{\alpha_\eta}+1\leq\gamma\leq\delta_{\alpha_{\eta+1}}$)}\\
\displaystyle\{\omega_1\xi\mid\xi\in C_{\alpha_\eta}\}\cup\bigcup_{\xi\in\mathrm{acc}(C_{\alpha_\eta})\cap\alpha_\eta}[\omega_1\xi, \omega_1\xi+\omega\delta_\xi)\\
\quad\text{(for $\alpha=\omega_1\alpha_\eta$, $\eta\in\mathrm{Lim}\cap\omega_1$).}
\end{cases}
$$
Then it is easy to check that $\seq{\overline{C}_\alpha\mid\alpha\in C}$ is end-extending.\qed

\section{Indestructible and absolute properness}\label{sec:doublestargames}

As we mentioned after Theorem \ref{thm:scpap}, $**$-tactically closed forcing does not necessarily preserve $\mathrm{PFA}$.
In this section we take a closer look on the extent of preservation of fragments of $\mathrm{PFA}$ under $**$-tactically closed forcing. For this purpose, we newly introduce two strengthenings of properness, to define corresponding fragments of $\mathrm{PFA}$.

\begin{dfn}
\begin{enumerate}[(1)]
  \item A poset $\mathbb{P}$ is said to be {\it indestructibly proper} if $\mathbb{P}\times\mathbb{Q}$ is proper for every $\sigma$-closed poset $\mathbb{Q}$.
  \item A poset $\mathbb{P}$ is said to be {\it absolutely proper} if $\mathbb{P}$ is proper in $V^\mathbb{Q}$ for every $\sigma$-closed poset $\mathbb{Q}$.
\end{enumerate}
\end{dfn}

\begin{lma}
\begin{enumerate}[(1)]
  \item\label{item:absind} Every absolutely proper poset is indestructibly proper.
  \item\label{item:baireind} Every totally proper ({\it i.e.} $\omega_1$-Baire and proper) poset is indestructibly proper.
  \item\label{item:scabs} If $\mathbb{P}=\mathbb{P}_0*\dot{\mathbb{P}}_1$ where $\mathbb{P}_0$ is $\sigma$-closed and $\dot{\mathbb{P}}_1$ is c.c.c. in $V^{\mathbb{P}_0}$ then $\mathbb{P}$ is absolutely proper.
\end{enumerate}
\end{lma}

\proof For (\ref{item:absind}), suppose $\mathbb{P}$ is absolutely proper and $\mathbb{Q}$ is $\sigma$-closed. Then $\mathbb{P}\times\mathbb{Q}\simeq\mathbb{Q}\times\mathbb{P}$ can be seen as a two step iteration of proper forcing and thus itself is proper.

For (\ref{item:baireind}), suppose $\mathbb{P}$ is totally proper and $\mathbb{Q}$ is $\sigma$-closed. Since $\mathbb{P}$ adds no new countable sequence of $\mathbb{Q}$-conditions, $\mathbb{Q}$ remains $\sigma$-closed in $V^{\mathbb{P}}$. Therefore $\mathbb{P}\times\mathbb{Q}$ can be seen as a two step iteration of proper forcing and thus itself proper.

For (\ref{item:scabs}), suppose that $\mathbb{P}=\mathbb{P}_0*\dot{\mathbb{P}}_1$ where $\mathbb{P}_0$ is $\sigma$-closed and $\dot{\mathbb{P}}_1$ is c.c.c. in $V^{\mathbb{P}_0}$ and that $\mathbb{Q}$ is $\sigma$-closed. Note first that $\mathbb{P}$ is always forcing equivalent to the two step iteration $\mathbb{P}_0*\dot{\mathbb{P}}_1$ in any extension of $V$. Thus to prove that $\mathbb{P}$ remains proper in $V^\mathbb{Q}$, it is enough to show that $\mathbb{P}_0$ is proper in $V^\mathbb{Q}$ and $\dot{\mathbb{P}}_1$ is proper in $(V^\mathbb{Q})^{\mathbb{P}_0}$. Since $\mathbb{Q}$ adds no new countable sequence of $\mathbb{P}_0$-conditions, $\mathbb{P}_0$ remains $\sigma$-closed in $V^{\mathbb{Q}}$. By the same argument we have $\mathbb{Q}$ is $\sigma$-closed in $V^{\mathbb{P}_0}$. By our assumption $\dot{\mathbb{P}}_1$ is c.c.c. in $V^{\mathbb{P}_0}$. Note that any c.c.c. poset remains c.c.c. under any $\sigma$-closed forcing, because any uncountable antichain in a generic extension by a $\sigma$-closed forcing can be approximated within the ground model. Applying this to $V^{\mathbb{P}_0}$ we have that $\dot{\mathbb{P}}_1$ remains c.c.c. in $(V^{\mathbb{P}_0})^\mathbb{Q}\simeq V^{\mathbb{P}_0\times\mathbb{Q}} \simeq (V^\mathbb{Q})^{\mathbb{P}_0}$. This completes our proof of (\ref{item:scabs}).\qed

It is known that there exists a poset which is indestructibly proper but not absolutely proper (see \cite{yoshinobu:_properness} for proof). We will take a closer look on the difference between these notions.

\begin{thm}\label{thm:maap}
Assume $\mathrm{PFA}$. Then for any $**$-tactically closed poset $\mathbb{P}$,
$$
\force_{\mathbb{P}}\mathrm{MA}_{\omega_1}(\text{absolutely proper})
$$
holds.\footnote{A weaker statement that $\mathbb{P}$ forces $\mathrm{MA}_{\omega_1}(\text{$\sigma$-closed*c.c.c.})$ under the same assumption has been announced in \cite{yoshinobu:_further} without proof.} In particular, $\mathrm{SCL^-_e}$ is consistent with $\mathrm{MA}_{\omega_1}(\text{absolutely proper})$.
\end{thm}

\proof Assume $\mathrm{PFA}$. Let $\mathbb{P}$ be any $**$-tactically closed poset, $\tau$ a winning $*$-tactic for $G^{**}(\mathbb{P})$, and $\dot{\mathbb{Q}}$ a $\mathbb{P}$-name for an absolutely proper poset in $V^\mathbb{P}$. Suppose $\dot{D}_\gamma$ is a $\mathbb{P}$-name for a dense subset of $\dot{\mathbb{Q}}$ for each $\gamma<\omega_1$. It is enough to show that for each $p_0\in\mathbb{P}$ there exists $p'\leq_{\mathbb{P}}p_0$ such that
$$
p'\Vdash_{\mathbb{P}}\text{\lq\lq$\exists F(\text{$F$ is  a filter on $\dot{\mathbb{Q}}$} \land\forall\gamma<\omega_1(F\cap\dot{D}_\gamma\not=\emptyset))$\rq\rq}.
$$
\begin{dfn}\label{dfn:tildep}
We define a poset $\tilde{\mathbb{P}}$ as follows: $\tilde{\mathbb{P}}$ consists of all conditions of the form $r=\seq{A^r_\gamma\mid\gamma\leq\beta^r}$ satisfying:
\begin{enumerate}[(1)]
\item\label{item:beta} $\beta^r<\omega_1$,
\item\label{item:ctble} $A^r_\gamma\in\mathbb{P}^*$ for each $\gamma\leq\beta^r$,
\item\label{item:Aprop} for each $\gamma<\delta\leq\beta^r$
  \begin{enumerate}[(i)]
    \item\label{subitem:inclusion} $A^r_\gamma\subseteq A^r_\delta$, 
    \item\label{subitem:doublestar} $p\leq_{\mathbb{P}}\tau(A^r_\gamma)$ for each $p\in A^r_\delta\setminus A^r_\gamma$ and
    \item\label{subitem:inequality} every common extension of $A^r_\delta$ also extends $\tau(A^r_\gamma)$.
  \end{enumerate}
\end{enumerate}
$\tilde{\mathbb{P}}$ is ordered by end-extension.
\end{dfn}

It almost appears that a $\tilde{\mathbb{P}}$-condition forms a sequence of Player $\mathrm{I}$'s moves in a play of $G^{**}(\mathbb{P})$ where Player $\mathrm{II}$ plays with $\tau$, but note that $\tilde{\mathbb{P}}$-conditions are not required to be $\subseteq$-continuous, unlike the rule of $G^{**}(\mathbb{P})$.

\begin{lma}\label{lma:ddensity}
\begin{enumerate}[(1)]
  \item\label{item:dextend} Every $\tilde{\mathbb{P}}$-condition can be properly extended.
  \item\label{item:dclosed} For any strictly decreasing sequence $\seq{r_n\mid n<\omega}$ in $\tilde{\mathbb{P}}$, where $r_n=\seq{A_\gamma\mid\gamma\leq\beta_n}$ for each $n<\omega$,
$$
\concat{(\bigcup_{n<\omega} r_n)}{\seq{\bigcup_{n<\omega}A_{\beta_n}}}
$$
is in $\tilde{\mathbb{P}}$ and extends all $r_n$'s. In particular $\tilde{\mathbb{P}}$ is $\sigma$-closed.
  \item\label{item:ddensity} For every $\beta<\omega_1$, $\tilde{E}_\beta=\{r\in\tilde{\mathbb{P}}\mid\beta^r\geq\beta\}$ is dense in $\tilde{\mathbb{P}}$.
\end{enumerate}
\end{lma}

\proof For (\ref{item:dextend}), let $r=\seq{A^r_\gamma\mid\gamma\leq\beta^r}\in\tilde{\mathbb{P}}$ be arbitrary. Then
$$
\concat{r}{\seq{A^r_{\beta^r}\cup\{\tau(A^r_{\beta^r})\}}}
$$
is in $\tilde{\mathbb{P}}$ and extends $r$. 

For (\ref{item:dclosed}), let $\seq{r_n\mid n<\omega}$ be any strictly decreasing sequence in $\tilde{\mathbb{P}}$. Set $r_n=\seq{A_\gamma\mid\gamma\leq\beta_n}$ for each $n<\omega$. Then
$$
\begin{matrix}
\text{Player}\ \rmi: & A_{\beta_0}\phantom{\tau(A_{\beta_0})} &  A_{\beta_1}\phantom{\tau(A_{\beta_1}} & A_{\beta_2}\phantom{\tau(A_{\beta_2}} & \cdots \\
\text{Player}\ \rmii: & \phantom{A_{\beta_0}}\tau(A_{\beta_0}) & \phantom{A_{\beta_1}}\tau(A_{\beta_1}) & \phantom{A_{\beta_2}}\tau(A_{\beta_2}) & \cdots
\end{matrix}
$$
forms a play of $G^{**}(\mathbb{P})$ where Player $\mathrm{II}$ plays with $\tau$, and thus $\bigcup_{n<\omega}A_{\beta_n}$ has a common extension in $\mathbb{P}$. Therefore
$$
\concat{(\bigcup_{n<\omega} r_n)}{\seq{\bigcup_{n<\omega}A_{\beta_n}}}
$$
is in $\tilde{\mathbb{P}}$ and extends all $r_n$'s.
 
(\ref{item:ddensity}) can be proved by a straightforward induction on $\beta$ using (\ref{item:dextend}) and (\ref{item:dclosed}).\qed

By Lemma \ref{lma:ddensity}, for any $\tilde{\mathbb{P}}$-generic filter $\tilde{G}$ over $V$, $\bigcup\tilde{G}$ is a sequence of length $\omega_1$, consisting of elements of $\mathbb{P}^*$.

\begin{dfn}
\begin{enumerate}[(1)]
\item For each $\gamma<\omega_1$, let $\dot{A}_\gamma$ be the $\tilde{\mathbb{P}}$-name representing the $\gamma$-th entry of $\bigcup\tilde{G}$ whenever $\tilde{G}$ is a $\tilde{\mathbb{P}}$-generic filter over $V$.
\item Let $\dot{S}$ be the $\tilde{\mathbb{P}}$-name such that
$$
  \force_{\tilde{\mathbb{P}}}\text{\lq\lq$\dot{S}=\{\gamma<\omega_1\mid \dot{A}_\gamma=\bigcup_{\xi<\gamma}\dot{A}_\xi\}$.\rq\rq}
$$
\end{enumerate}
\end{dfn}

\begin{lma}\label{lma:stat}
$\dot{S}$ is a stationary subset of $\omega_1$ in $V^{\tilde{\mathbb{P}}}$.
\end{lma}

\proof
%Note that $\mathbb{P}*\dot{\mathbb{Q}}*\dot{\mathbb{R}}$ is forcing equivalent to $\mathbb{P}*\dot{\mathbb{R}}*\dot{\mathbb{Q}}$. Since $\dot{\mathbb{Q}}$ is absolutely proper in $V^\mathbb{P}$ and $\dot{\mathbb{R}}$ is $\sigma$-closed in $V^\mathbb{P}$ by Lemme \ref{lma:rdensity}(\ref{item:closed}), $\dot{\mathbb{Q}}$ is proper in $V^{\mathbb{P}*\dot{\mathbb{R}}}$. Therefore it is enough to show that $\dot{S}$ is a stationary subset of $\omega_1$ in $V^{\mathbb{P}*\dot{\mathbb{R}}}$.
Let $\dot{E}$ be any $\tilde{\mathbb{P}}$-name for a club subset of $\omega_1$, and $r\in\tilde{\mathbb{P}}$ be arbitrary. Let $\theta$ be a sufficiently large regular cardinal, and $N$ a countable elementary submodel of $H_\theta$ containing all relevant objects. Let $\delta=N\cap\omega_1$. Pick an $(N, \tilde{\mathbb{P}})$-generic sequence $\seq{r_n\mid n<\omega}$ below $r$, and let $r_n=\seq{A_\gamma\mid\gamma\leq\beta_n}$. By genericity and Lemma \ref{lma:ddensity}(\ref{item:ddensity}) we have $\sup_{n<\omega}\beta_n=\delta$. By Lemma \ref{lma:ddensity}(\ref{item:dclosed}), $r_\omega=\concat{\bigcup_{n<\omega}r_n}{\seq{\bigcup_{n<\omega}A_{\beta_n}}}$ is a $\tilde{\mathbb{P}}$-condition extending $r$ and is $(N, \tilde{\mathbb{P}})$-generic. Now it is easy to see that $r_\omega$ forces that $\delta\in\dot{E}\cap\dot{S}\not=\emptyset$.\qed

\begin{dfn}\label{dfn:r}
\begin{enumerate}[(1)]
\item Let $\pi:\tilde{\mathbb{P}}\to\mathbb{P}$ be defined as follows: for each $r=\seq{A^r_\gamma\mid\gamma\leq\beta^r}\in\tilde{\mathbb{P}}$ we set
$$
\pi(r)=\tau(A^r_{\beta^r}).
$$
\item Let $\dot{\mathbb{R}}$ be the $\mathbb{P}$-name such that
$$
\dot{\mathbb{R}}_G=\{r\in\tilde{\mathbb{P}}\mid\pi(r)\in G\}\ \text{(ordered by $\leq_{\tilde{\mathbb{P}}}$)}
$$
for any $\mathbb{P}$-generic $G$ over $V$.
\end{enumerate}
\end{dfn}

It is easily checked that $\pi$ is a projection in the sense of \cite{cummings:_iterated}, and by the general theory of projection $\tilde{\mathbb{P}}$ is densely embedded into $\mathbb{P}*\dot{\mathbb{R}}$ by the map $\iota:\tilde{\mathbb{P}}\to\mathbb{P}*\dot{\mathbb{R}}$ defined by
$$
\iota(r)=\seq{\pi(r), \check{r}}.
$$

\begin{lma}\label{lma:r}
$\dot{\mathbb{R}}$ is $\sigma$-closed in $V^{\mathbb{P}}$.
\end{lma}
\proof Let $G$ be any $\mathbb{P}$-generic filter over $V$. In $V[G]$, suppose $\seq{r_n\mid n<\omega}$ is a strictly decreasing sequence in $\dot{\mathbb{R}}_G$. We may write as $r_n=\seq{A_\gamma\mid\gamma\leq\beta_n}$, where $\seq{\beta_n\mid n<\omega}$ is strictly increasing. Now let $\beta=\sup_{n<\omega}\beta_n$ and $A=\bigcup_{\gamma<\beta}A_\gamma$. Note that $A\subseteq G$ holds, and that $A\in V$ holds since $\mathbb{P}$ is $\omega_1$-Baire. Therefore $A$ has a common extension. Since the set
$$\{\tau(A\cup\{p\})\mid\text{$p$ is a common extension of $A$}\}$$
is in $V$ and is dense in the set of common extensions of $A$, we can choose $p$ such that $p$ is a common extension of $A$ and $\tau(A\cup\{p\})\in G$. Then by letting $A_\beta=A\cup\{p\}$ we have that $\seq{A_\gamma\mid\gamma\leq\beta}$ is an $\dot{\mathbb{R}}_G$-condition and extends all $r_n$'s.
\qed

We can naturally regard as $V^{\mathbb{P}}\subseteq V^{\mathbb{P}*\dot{\mathbb{R}}}=V^{\tilde{\mathbb{P}}}$. Since $\dot{\mathbb{Q}}$ is absolutely proper in $V^{\mathbb{P}}$ and $\dot{\mathbb{R}}$ is $\sigma$-closed in $V^{\mathbb{P}}$ by Lemma \ref{lma:r}, $\dot{\mathbb{Q}}$ remains proper in $V^{\tilde{\mathbb{P}}}$. By Lemma \ref{lma:stat} $\dot{S}$ is stationary in $V^{\tilde{\mathbb{P}}}$, and since $\dot{\mathbb{Q}}$ is proper in $V^{\tilde{\mathbb{P}}}$, $\dot{S}$ remains stationary in $V^{\tilde{\mathbb{P}}*\dot{\mathbb{Q}}}$.

\begin{dfn}
Let $\dot{\mathbb{S}}$ be a $(\tilde{\mathbb{P}}*\dot{\mathbb{Q}})$-name for the poset consisting of all closed bounded subsets of $\dot{S}$ and ordered by end-extension.
\end{dfn}

As it is well known, $\dot{\mathbb{S}}$ is $\sigma$-Baire in $V^{\tilde{\mathbb{P}}*\dot{\mathbb{Q}}}$ and forcing with $\dot{\mathbb{S}}$ adds a club subset of $\dot{S}$. Let $\dot{C}$ denote the $(\tilde{\mathbb{P}}*\dot{\mathbb{Q}}*\dot{\mathbb{S}})$-name for this club subset.

\begin{lma}
$\tilde{\mathbb{P}}*\dot{\mathbb{Q}}*\dot{\mathbb{S}}$ is proper.
\end{lma}

\proof Let $\seq{r, \dot{q}, \dot{s}}\in\tilde{\mathbb{P}}*\dot{\mathbb{Q}}*\dot{\mathbb{S}}$ be arbitrary. Let $\theta$ be a sufficiently large regular cardinal and $N$ a countable elementary submodel of $H_\theta$ containing all relevant objects. Let $\delta=N\cap\omega_1$. By an argument similar to the proof of Lemma  \ref{lma:stat} we have that there is $r'\leq_{\tilde{\mathbb{P}}}r$ which is $(N, \tilde{\mathbb{P}})$-generic and forces that $\delta\in\dot{S}$. For any $\tilde{\mathbb{P}}$-generic filter $\tilde{G}$ over $V$ with $r'\in\tilde{G}$, find an $(N[\tilde{G}], \dot{\mathbb{Q}}_{\tilde{G}})$-generic condition $q'\leq\dot{q}_{\tilde{G}}$. This is possible since $\dot{\mathbb{Q}}_{\tilde{G}}$ is proper in $V[\tilde{G}]$. Now for any $\dot{\mathbb{Q}}_{\tilde{G}}$-generic filter $H$ over $V[\tilde{G}]$ with $q'\in H$, pick an $(N[\tilde{G}][H], \dot{\mathbb{S}}_{\tilde{G}*H})$-generic sequence $\seq{s_n\mid n<\omega}$ extending $\dot{s}_{\tilde{G}*H}$. Since $\tilde{\mathbb{P}}*\dot{\mathbb{Q}}$ is proper, $N[\tilde{G}][H]\cap\omega_1=\delta$ and by genericity $\sup_{n<\omega}\max{s_n}=\delta$. Since $\delta\in\dot{S}_{\tilde{G}}$, $s'=\bigcup_{n<\omega}s_n\cup\{\delta\}\in\dot{\mathbb{S}}_{\tilde{G}*H}$ is an $(N[\tilde{G}][H], \dot{\mathbb{S}}_{\tilde{G}*H})$-generic condition extending $\dot{s}_{\tilde{G}*H}$. This completes the proof.\qed

We will denote $\tilde{\mathbb{P}}*\dot{\mathbb{Q}}*\dot{\mathbb{S}}$ as $\overline{\mathbb{P}}$. Pick any condition $\overline{p}_0=\seq{r_0, \dot{q}_0, \dot{s}_0}\in\overline{\mathbb{P}}$ such that $p_0\in A^{r_0}_0$ (remember that $p_0$ was given in the beginning of this proof). Now let
$$
D=\{\seq{r, \dot{q}, \dot{s}}\in\overline{\mathbb{P}}\mid\text{$\dot{q}$ is a $\mathbb{P}$-name}\}.
$$
Note that $D$ is a dense subset of $\overline{\mathbb{P}}$. Moreover for each $\gamma<\omega_1$ let
\begin{eqnarray*}
D_\gamma&=&\{\seq{r, \dot{q}, \dot{s}}\in D\mid r\force_{\tilde{\mathbb{P}}}\text{\lq\lq$\dot{q}\in\dot{D}_\gamma$\rq\rq}\},\\
E_\gamma&=&\{\seq{r, \dot{q}, \dot{s}}\in\overline{\mathbb{P}}\mid\beta^r\geq\gamma\},\\
U_\gamma&=&\{\overline{p}\in\overline{\mathbb{P}}\mid\exists\eta(\gamma\leq\eta<\omega_1\land\overline{p}\force_{\overline{\mathbb{P}}}\text{\lq\lq$\eta\in\dot{C}$\rq\rq})\}\ \text{and}\\
C_\gamma&=&\{\overline{p}\in\overline{\mathbb{P}}\mid\overline{p}\force_{\overline{\mathbb{P}}}\text{\lq\lq$\omega\gamma\in\dot{C}$\rq\rq}\lor\exists\xi<\omega\gamma(\overline{p}\force_{\overline{\mathbb{P}}}\text{\lq\lq$[\xi, \omega\gamma)\cap\dot{C}=\emptyset$\rq\rq})\}.
\end{eqnarray*}

It is straightforward to check that $D_\gamma$, $E_\gamma$, $U_\gamma$ and $C_\gamma$ for limit $\gamma$ are all dense subsets of $\overline{\mathbb{P}}$. Now by $\mathrm{PFA}$ there exists a filter $\mathcal{F}$ on $\overline{\mathbb{P}}$ satisfying
\begin{enumerate}[(1)]
  \item $\overline{p}_0\in\mathcal{F}$,
  \item $\mathcal{F}$ intersects with $D_\gamma$, $E_\gamma$, $U_\gamma$ and $C_\gamma$ for every $\gamma<\omega_1$, and 
  \item $\mathcal{F}$ is generated by $\mathcal{F}\cap D$.
\end{enumerate}
  Since conditions in $\mathcal{F}$ are pariwise compatible and $\mathcal{F}$ intersects with all $E_\gamma$'s, we may write as
$$
\bigcup\{r\mid\seq{r, \dot{q}, \dot{s}}\in\mathcal{F}\ \text{for some $\dot{q}$ and $\dot{s}$}\}=\seq{A_\gamma\mid\gamma<\omega_1}.
$$
Note that each proper initial segment of $\seq{A_\gamma\mid\gamma<\omega_1}$ is decided by a single condition in $\mathcal{F}$ and therefore satisfies
\begin{enumerate}[(1)]
\item\label{item:fctble} $A_\gamma\in\mathbb{P}^*$ for each $\gamma<\omega_1$,
\item\label{item:fAprop} for each $\gamma<\delta<\omega_1$
  \begin{enumerate}[(i)]
    \item\label{subitem:finclusion} $A_\gamma\subseteq A_\delta$, 
    \item\label{subitem:fdoublestar} $p\leq_{\mathbb{P}}\tau(A_\gamma)$ for each $p\in A_\delta\setminus A_\gamma$,
    \item\label{subitem:finequality} every common extension of $A_\delta$ also extends $\tau(A_\gamma)$.
  \end{enumerate}
\end{enumerate}

Now set
$$
C=\{\gamma<\omega_1\mid\exists\overline{p}\in\mathcal{F}(\overline{p}\force_{\overline{\mathbb{P}}}\text{\lq\lq$\gamma\in\dot{C}$\rq\rq})\}.
$$
Note that $C$ is unbounded in $\omega_1$ since $\mathcal{F}$ intersects with all $U_\gamma$'s, and is closed since $\mathcal{F}$ intersects with all $C_\gamma$'s. Let $\seq{\gamma_\zeta\mid\zeta<\omega_1}$ be the increasing enumeration of $C$. For each $\gamma\in C$, in fact there exists $\overline{p}\in\mathcal{F}\cap E_\gamma$ with $\overline{p}\force_{\overline{\mathbb{P}}}\text{\lq\lq$\gamma\in\dot{C}$\rq\rq}$, and thus $A_\gamma=\bigcup_{\xi<\gamma}A_\xi$ holds. In particular for each limit $\zeta<\omega_1$ we have
$$
A_{\gamma_\zeta}=\bigcup_{\xi<\zeta}A_{\gamma_\xi}
$$
and thus
$$
\begin{matrix}
\text{Player}\ \rmi: & A_{\gamma_0}\phantom{\tau(A_{\gamma_0})} &  A_{\gamma_1}\phantom{\tau(A_{\gamma_1})} & \cdots & A_{\gamma_\omega}\phantom{\tau(A_{\gamma_\omega})} & A_{\gamma_{\omega+1}}\phantom{\tau(A_{\gamma_{\omega+1}})} & \cdots \\
\text{Player}\ \rmii: & \phantom{A_{\gamma_0}}\tau(A_{\gamma_0}) & \phantom{A_{\gamma_1}}\tau(A_{\gamma_1}) & \cdots & \phantom{A_\omega}\tau(A_{\gamma_\omega}) &  \phantom{A_{\omega+1}}\tau(A_{\gamma_{\omega+1}}) & \cdots
\end{matrix}
$$
forms a play of $G^{**}(\mathbb{P})$ where Player $\mathrm{II}$ plays with $\tau$. Therefore $\bigcup_{\xi<\omega_1}A_\xi$ has a common extension $p'\in\mathbb{P}$. Clearly $p'\leq_{\mathbb{P}}p_0$ holds. Now since $\mathcal{F}\cap D$ generates $\mathcal{F}$ and thus is directed, for every $\mathbb{P}$-generic filter $G$ over $V$ with $p'\in G$,
$$
\mathcal{D}=\{\dot{q}_G\mid\seq{r, \dot{q}, \dot{s}}\in\mathcal{F}\cap D\ \text{for some $r$ and $\dot{s}$}\}
$$
is a directed subset of $\dot{\mathbb{Q}}_G$ intersecting with $(\dot{D}_\gamma)_G$ for every $\gamma<\omega_1$, and thus generates a filter $F$ as desired. This completes the proof of Theorem \ref{thm:maap}.\qed

Note that even $\mathrm{MA}_{\omega_1}(\text{$\sigma$-closed*c.c.c.})$ implies many consequences of $\mathrm{PFA}$. For example, it negates $\square_\kappa$ for every uncountable cardinal $\kappa$ (see \cite[Theorem 3.9 and Corollary 3.11]{bekkali} for proof).\footnote{The failure of $\square_\kappa$ for every uncountable cardinal $\kappa$ under $\mathrm{PFA}$ was originally proved by Todorcevic \cite{PFAnote}.} It also implies the tree property of $\omega_2$ (due to Baumgartner; see \cite{devlin:_yorkshireman} for proof). Theorem \ref{thm:maap} thus shows that $\mathrm{SCL}^-_e$ is consistent with a considerably large fragment of $\mathrm{PFA}$, though is is not consistent with full $\mathrm{PFA}$.

As a sharpening of this contrast we have
\begin{thm}\label{thm:ind}
$\mathrm{MA}_{\omega_1}(\text{indestructibly proper})$ negates $\mathrm{SCL^-_e}$.
\end{thm}

To prove Theorem \ref{thm:ind}, we use the Mapping Reflection Principle ($\mathrm{MRP}$) introduced by Moore.

\begin{dfn}[Moore \cite{moore05:_set}]
\begin{enumerate}[(1)]
\item For a countable set $x$ and $a\in[x]^{<\omega}$, $[a, x]$ denotes the set $\{y\in[x]^\omega\mid a\subseteq y\}$.
\item For a regular cardinal $\theta\geq\omega_2$, an uncountable set $X\in H_\theta$ and a club subset $\mathcal{E}$ of $[H_\theta]^\omega$, a function $\Sigma:\mathcal{E}\to\mathcal{P}([X]^\omega)$ is said to be an {\it open stationary set mapping\/} if for every $M\in\mathcal{E}$,
\begin{enumerate}[(i)]
  \item $\Sigma(M)$ is {\it Ellentuck open\/}, that is, for every $x\in\Sigma(M)$ there is $a\in[x]^{<\omega}$ such that $[a, x]\subseteq\Sigma(M)$, and
  \item $\Sigma(M)$ is {\it $M$-stationary\/}, that is, $C\cap M\cap\Sigma(M)\not=\emptyset$ holds for every club subset $C$ of $[X]^\omega$ in $M$.
\end{enumerate}
\item For an open stationary set mapping $\Sigma:\mathcal{E}\to\mathcal{P}([X]^\omega)$, a sequence $\seq{N_\xi\mid\xi<\omega_1}$ is said to be a {\it reflecting sequence\/} for $\Sigma$ if
\begin{enumerate}[(i)]
  \item for every $\xi<\omega_1$, $N_\xi\in\mathcal{E}$ and $N_\xi\in N_{\xi+1}$ holds, and
  \item for every limit $\xi<\omega_1$, $N_\xi=\bigcup_{\nu<\xi}N_\nu$ holds and there exists $\eta<\xi$ such that $N_\nu\cap X\in\Sigma(N_\xi)$ holds for every $\nu\in(\eta, \xi)$.
\end{enumerate}
  \item The {\it Mapping Reflection Principle\/} ($\mathrm{MRP}$) is the statement that every open stationary set mapping has a reflecting sequence for it.
\end{enumerate}
\end{dfn}

From the proof that $\mathrm{PFA}$ implies $\mathrm{MRP}$ in \cite{moore05:_set}, we can observe the following theorem.
\begin{thm}\label{thm:tot}
$\mathrm{MA}_{\omega_1}(\text{totally proper})$ implies $\mathrm{MRP}$ (and thus so does $\mathrm{MA}_{\omega_1}(\text{indestructibly proper})$).
\end{thm}

By Theorem \ref{thm:tot}, for Theorem \ref{thm:ind} it is enough to show the following.

\begin{lma}\label{lma:mrp}
$\mathrm{MRP}$ negates $\mathrm{SCL^-_e}$.
\end{lma}

\proof Assume $\mathrm{MRP}$. We will show that $\mathrm{SCL^-_e}$ fails. Suppose to the contrary that there exists an $\mathrm{SCL^-_e}$-system $\vec{C}=\seq{C_\alpha\mid \alpha\in S^2_0}$. Let $\theta$ be a sufficiently large regular cardinal. Fix a complete set of Skolem functions of the structure $\mathcal{A}=\seq{H_\theta, \in, \vec{C}}$ so that we can define the Skolem hull operation $\mathrm{Sk}^\mathcal{A}$. For each $\alpha\in S^2_0$, let
$$
t_\alpha=\{\eta\in S^2_0\cap\alpha\mid C_\eta=C_\alpha\cap\eta\}
$$
and
$$
\mathcal{F}_\alpha=\{\mathrm{Sk}^\mathcal{A}(t_\alpha\cap\beta)\cap\omega_2\mid\beta<\alpha\}.
$$
Note that each $\mathcal{F}_\alpha$ is a subset of $[\omega_2]^\omega$ and is linearly ordered by $\subseteq$. Now let
$$
\mathcal{E}=\{M\in[H_\theta]^\omega\mid\mathrm{Sk}^\mathcal{A}(M)=M\}.
$$
We will define an open stationary set mapping $\Sigma:\mathcal{E}\to\mathcal{P}([\omega_2]^\omega)$. Pick any $M\in\mathcal{E}$. Let $\alpha_M=\sup(M\cap\omega_2)$ and let $\seq{D^M_n\mid n<\omega}$ enumerate the club subsets of $[\omega_2]^\omega$ in $M$. For each $n<\omega$, since $D^M_n$ is a club subset of $[\omega_2]^\omega$, there exist $x^{M, n}_0$, $x^{M, n}_1\in D^M_n$ which are $\subseteq$-incomparable. By elementarity we can choose such $x^{M, n}_0$, $x^{M, n}_1$ in $M$. Now pick $\xi^{M, n}_i\in x^{M, n}_i\setminus x^{M, n}_{1-i}$ and let $I^{M, n}_i=[\{\xi^{M, n}_i\}, x^{M, n}_i]$ for each $i<2$. Then note that any pair of members respectively of $I^{M, n}_0$ and $I^{M, n}_1$ are $\subseteq$-incomparable. Since $\mathcal{F}_{\alpha_M}$ is linearly ordered by $\subseteq$, at least one of $I^{M, n}_0$ and $I^{M, n}_1$ does not intersect with $\mathcal{F}_{\alpha_M}$. Let $i^{M, n}<2$ be such that $I^{M, n}_{i^{M, n}}\cap\mathcal{F}_{\alpha_M}=\emptyset$. Now set
$$
\Sigma(M)=\bigcup_{n<\omega}I^{M, n}_{i^{M, n}}.
$$
Clearly $\Sigma(M)$ is Ellentuck open and $M$-stationary. Moreover it holds that
$$
\Sigma(M)\cap\mathcal{F}_{\alpha_M}=\emptyset.
$$

Now apply $\mathrm{MRP}$ and let $\vec{N}=\seq{N_\xi\mid\xi<\omega_1}$ be a reflecting sequence for $\mathcal{E}$. Note that $\vec{N}$ forms an elementary chain of countable elementary submodels of $\mathcal{A}$, and thus $N=\bigcup_{\xi<\omega_1}N_\xi$ is also an elementary submodel of $\mathcal{A}$. Moreover since both $\seq{N_\xi\cap\omega_1\mid\xi<\omega_1}$ and $\seq{\alpha_{N_\xi}\mid\xi<\omega_1}$ are strictly increasing, it holds that $\omega_1\subseteq N$ and that $\gamma=N\cap\omega_2\in S^2_1$. Now let $C$ be a club subset of $\gamma\cap S^2_0$ with $\mathrm{o.t.}(C)=\omega_1$ such that $\seq{C_\alpha\mid\alpha\in C}$ is end-extending. Now let
$$
t=\{\nu\in\gamma\cap S^2_0\mid(\bigcup_{\alpha\in C}C_\alpha)\cap\nu=C_\nu\}.
$$
Note that $C\subseteq t\subseteq\gamma$ holds. Note also that for each $\alpha\in t$, it holds that $t\cap\alpha=t_\alpha$ and thus $\mathrm{o.t.}(t\cap\alpha)$ is definable from $\alpha$ in $\mathcal{A}$. This shows that $\omega_1\subseteq\mathrm{Sk}^\mathcal{A}(t)$ and thus $\mathrm{Sk}^\mathcal{A}(t)\cap\omega_2$ is an ordinal. By the definition of $\gamma$, it holds that $\mathrm{Sk}^\mathcal{A}(t)\cap\omega_2\subseteq\mathrm{Sk}^\mathcal{A}(\gamma)\cap\omega_2=\gamma$, and since $t$ is unbounded in $\gamma$ we have $\mathrm{Sk}^\mathcal{A}(t)\cap\omega_2=\gamma$. Now let
$$
\mathcal{F}=\{\mathrm{Sk}^\mathcal{A}(t\cap\alpha)\cap\omega_2\mid\alpha\in C\}.
$$
Then both $\mathcal{F}$ and $\seq{N_\xi\cap\omega_2\mid\xi<\omega_1}$ are club subsets of $[\gamma]^\omega$, and thus there exists a club subset $D$ of $\omega_1$ such that $N_\nu\cap\omega_2\in\mathcal{F}$ for every $\nu\in D$. Since $\vec{N}$ is a reflecting sequence, by Fodor's Lemma there exists $\eta_0<\omega_1$ and a stationary subset $S$ of $\omega_1\cap\mathrm{Lim}$ such that $N_\nu\cap\omega_2\in\Sigma(N_\xi)$ for every $\nu\in(\eta_0, \omega_1)$ and $\xi\in S\setminus(\nu+1)$. Now pick $\nu\in D\setminus(\eta_0+1)$. Then on the one hand $N_\nu\cap\omega_2\in\mathcal{F}$ holds by the choice of $D$, but on the other hand $N_\nu\cap\omega_2\notin\mathcal{F}_{\alpha_{N_\xi}}$ holds for every $\xi\in S\setminus(\nu+1)$ and therefore it holds that $N_\nu\cap\omega_2\notin\mathcal{F}$. This is a contradiction.\qed

The following is another result contrasting with Theorem \ref{thm:ind}.

\begin{thm}\label{thm:maip}
Assume $\mathrm{PFA}$. Then for any $(\omega_1+1)$-strongly strategically closed poset $\mathbb{P}$,
$$
\force_{\mathbb{P}}\mathrm{MA}_{\omega_1}(\text{indestructibly proper})
$$
holds.
\end{thm}

\proof The proof is to a certain extent similar to that of Theorem \ref{thm:maap}, but simpler than that. Assume $\mathrm{PFA}$. Let $\mathbb{P}$ be any $(\omega_1+1)$-strongly strategically closed poset, $\sigma$ a winning strategy for $G^{\mathrm{I}}_{\omega_1+1}(\mathbb{P})$, and $\dot{\mathbb{Q}}$ a $\mathbb{P}$-name for an indestructibly proper poset in $V^\mathbb{P}$. Suppose $\dot{D}_\gamma$ is a $\mathbb{P}$-name for a dense subset of $\dot{\mathbb{Q}}$ for each $\gamma<\omega_1$. Again it is enough to show that for each $p_0\in\mathbb{P}$ there exists $p'\leq_{\mathbb{P}}p_0$ such that
$$
p'\Vdash_{\mathbb{P}}\text{\lq\lq$\exists F(\text{$F$ is  a filter on $\dot{\mathbb{Q}}$} \land\forall\gamma<\omega_1(F\cap\dot{D}_\gamma\not=\emptyset))$\rq\rq}.
$$

\begin{dfn}\label{dfn:strongtildep}
We define a poset $\tilde{\mathbb{P}}$ as follows: $\tilde{\mathbb{P}}$ consists of all $\leq_\mathbb{P}$-descending sequence of the form $r=\seq{a^r_\gamma\mid\gamma\leq\beta^r}$ satisfying $\beta^r<\omega_1$ and $a^r_{\gamma+1}\leq_\mathbb{P}\sigma(r\upharpoonright(\gamma+1))$ for every $\gamma<\beta^r$. $\tilde{\mathbb{P}}$ is ordered by end-extension.
\end{dfn}

Note that a $\tilde{\mathbb{P}}$-condition forms a sequence of Player $\mathrm{I}$'s moves in a play of $G^{\mathrm{I}}_{\omega_1+1}(\mathbb{P})$ where Player $\mathrm{II}$ plays with $\sigma$ really this time. This easily implies facts about $\tilde{\mathbb{P}}$ similar to Lemma \ref{lma:ddensity}.

%: For every $r\in\tilde{\mathbb{P}}$, $\concat{r}{\seq{\sigma(r)}}$ is a $\tilde{\mathbb{P}}$-condition extending $r$, and  every descending sequence $\seq{r_n\mid n<\omega}$ in $\tilde{\mathbb{P}}$ has a common extension. Thus we have the following.

\begin{lma}\label{lma:strongddensity}
\begin{enumerate}[(1)]
  \item\label{item:strongdextend} Every $\tilde{\mathbb{P}}$-condition can be properly extended.
  \item\label{item:strongdclosed} $\tilde{\mathbb{P}}$ is $\sigma$-closed.
  \item\label{item:strongddensity} For every $\beta<\omega_1$, $\tilde{E}_\beta=\{r\in\tilde{\mathbb{P}}\mid\beta^r\geq\beta\}$ is dense in $\tilde{\mathbb{P}}$. \qed
\end{enumerate}
\end{lma}

\begin{dfn}\label{dfn:strongr}
Let $\dot{\mathbb{R}}$ be the $\mathbb{P}$-name such that
$$
\dot{\mathbb{R}}_G=\{r\in\tilde{\mathbb{P}}\mid\sigma(r)\in G\}\ \text{(ordered by $\leq_{\tilde{\mathbb{P}}}$)}
$$
for any $\mathbb{P}$-generic $G$ over $V$.
\end{dfn}

Again it is easy to see that $\sigma$ is a projection from $\tilde{\mathbb{P}}$ to $\mathbb{P}$ and that $\tilde{\mathbb{P}}$ is densely embedded into $\mathbb{P}*\dot{\mathbb{R}}$ by $\iota:\tilde{\mathbb{P}}\to\mathbb{P}*\dot{\mathbb{R}}$ defined by $\iota(r)=\seq{\sigma(r), \check{r}}$. We also have the following.

\begin{lma}\label{lma:strongr}
$\dot{\mathbb{R}}$ is $\sigma$-closed in $V^{\mathbb{P}}$.
\end{lma}

\proof Let $G$ be any $\mathbb{P}$-generic filter over $V$. In $V[G]$, suppose $\seq{r_n\mid n<\omega}$ is a strictly decreasing sequence in $\dot{\mathbb{R}}_G$. Since $\mathbb{P}$ is $\omega_1$-Baire, we have $\seq{r_n\mid n<\omega}\in V$. Therefore $r=\bigcup_{n<\omega} r_n$ lies in $V$ and forms a sequence of Player $\mathrm{I}$'s moves in a play of $G^{\mathrm{I}}_{\omega_1+1}(\mathbb{P})$ where Player $\mathrm{II}$ plays with $\sigma$, and thus has a common extension in $\mathbb{P}$. Therefore
$$
\{\sigma(\concat{r}{\seq{p}})\mid\text{$p$ is a common extension of $r$ in $\mathbb{P}$}\}
$$
is dense in the set of common extensions of $r$ in $\mathbb{P}$, and thus we can pick $p\in G$ which is a common extension of $r$ in $\mathbb{P}$. For this $p$ we have that $\concat{r}{\seq{p}}\in\dot{\mathbb{R}}_G$ and is a common extension of $\seq{r_n\mid n<\omega}$ in $\dot{\mathbb{R}}_G$. \qed

Since $\dot{\mathbb{Q}}$ is indestructibly proper in $V^\mathbb{P}$, by Lemma \ref{lma:strongr} $\dot{\mathbb{R}}\times\dot{\mathbb{Q}}$ is proper in $V^\mathbb{P}$. Since $\mathbb{P}$ itself is proper, so is $\mathbb{P}*(\dot{\mathbb{R}}\times\dot{\mathbb{Q}})$, which is forcing equivalent to $\mathbb{P}*\dot{\mathbb{R}}*\dot{\mathbb{Q}}=\tilde{\mathbb{P}}*\dot{\mathbb{Q}}$. Pick $\overline{p}_0=\seq{r_0, \dot{q}_0}\in\tilde{\mathbb{P}}*\dot{\mathbb{Q}}$ such that $a^{r_0}_0=p_0$. Now for each $\gamma<\omega_1$ let
\begin{eqnarray*}
D&=&\{\seq{r, \dot{q}}\in\tilde{\mathbb{P}}*\dot{\mathbb{Q}}\mid\text{$\dot{q}$ is a $\mathbb{P}$-name}\},\\
D_\gamma&=&\{\seq{r, \dot{q}}\in D\mid r\force_{\tilde{\mathbb{P}}}\text{\lq\lq$\dot{q}\in\dot{D}_\gamma$\rq\rq}\}\quad\text{and}\\
E_\gamma&=&\{\seq{r, \dot{q}}\in\tilde{\mathbb{P}}*\dot{\mathbb{Q}}\mid \beta^r\geq\gamma\}.
\end{eqnarray*}
Note that $D$, $D_\gamma$'s and $E_\gamma$'s are all dense subsets of $\tilde{\mathbb{P}}*\dot{\mathbb{Q}}$.
Now by $\mathrm{PFA}$ there exists a filter $\mathcal{F}$ on $\tilde{\mathbb{P}}*\dot{\mathbb{Q}}$ satisfying
\begin{enumerate}[(1)]
  \item $\overline{p}_0\in\mathcal{F}$,
  \item $\mathcal{F}$ intersects with $D_\gamma$ and $E_\gamma$ for every $\gamma<\omega_1$, and
  \item $\mathcal{F}$ is generated by $\mathcal{F}\cap D$.
\end{enumerate}
Exactly as in the proof of Theorem \ref{thm:maap}, we may write as
$$
\bigcup\{r\mid\seq{r, \dot{q}}\in\mathcal{F}\ \text{for some $\dot{q}$}\}=\seq{a_\gamma\mid\gamma<\omega_1}.
$$
Moreover we have that  $\seq{a_\gamma\mid\gamma<\omega_1}$ forms a sequence of Player $\rmi$'s moves in a play of $G^\rmi_{\omega_1+1}(\mathbb{P})$ where Player $\rmii$ plays with $\sigma$, and therefore has a common extension $p'\in\mathbb{P}$. Clearly $p'\leq_\mathbb{P}p_0$ holds. The rest of the proof is exactly as that of Theorem \ref{thm:maap}: For every $\mathbb{P}$-generic filter $G$ over $V$ with $p'\in G$, 
$$
\mathcal{D}=\{\dot{q}_G\mid\seq{r, \dot{q}}\in\mathcal{F}\cap D\ \text{for some $r$}\}
$$
is a directed subset of $\dot{\mathbb{Q}}_G$ which intersects with $(\dot{D_\gamma})_G$ for each $\gamma<\omega_1$, and thus generates a desired filter.\qed
%$$
%\mathbb{P}*\dot{\mathbb{Q}}*\dot{\mathbb{R}}\simeq\mathbb{P}*(\dot{\mathbb{Q}}\times\dot{\mathbb{R}})\simeq\mathbb{P}*(\dot{\mathbb{R}}\times\dot{\mathbb{Q}})\simeq\mathbb{P}*\dot{\mathbb{R}}*\dot{\mathbb{Q}}\simeq\tilde{\mathbb{P}}*\dot{\mathbb{Q}}.
%$$ 

\begin{cor}
$\mathrm{SCL}+\mathrm{AP}_{\omega_1}$ does not imply $\mathrm{SCL^-_e}$.
\end{cor}
\proof Start with $\mathrm{PFA}+\mathrm{SCL}$. Force $\mathrm{AP}_{\omega_1}$ with the poset $\mathbb{P}_{\mathrm{AP}_{\omega_1}}$. Since $\mathbb{P}_{\mathrm{AP}_{\omega_1}}$ is $(\omega_1+1)$-strongly strategically closed, by Theorem \ref{thm:maip} $\mathrm{MA}_{\omega_1}(\text{indestructibly proper})$ holds in the generic extension. By Theorem \ref{thm:ind}, it implies that $\mathrm{SCL^-_e}$ fails in the model. On the other hand, in the same model $\mathrm{AP}_{\omega_1}$ holds and $\mathrm{SCL}$ is preserved since $\mathbb{P}_{\mathrm{AP}_{\omega_1}}$ is $\omega_2$-Baire.\qed

%%%%%%%%%%%%%%%%%%%%%%%%%%%%%%%%%%%%%%%%%%%%
\section{Summary}\label{sec:summary}
%%%%%%%%%%%%%%%%%%%%%%%%%%%%%%%%%%%%%%%%%%%%
Implications between fragments of $\square_{\omega_1}$ mentioned in this paper are indicated in the left half of the following diagram. The variations of setwise climbability properties, which are newly introduced and argued in this paper, are indicated with rounded rectangles. Relevant properties of posets are listed in the right half of the diagram, and each broken line indicates that the combinatorial property at its left end is (naturally) equivalent to $\mathrm{MA}_{\omega_2}$ for the class of posets with the property at its right end.

\begin{figure}[h]
\includegraphics[height=9cm]{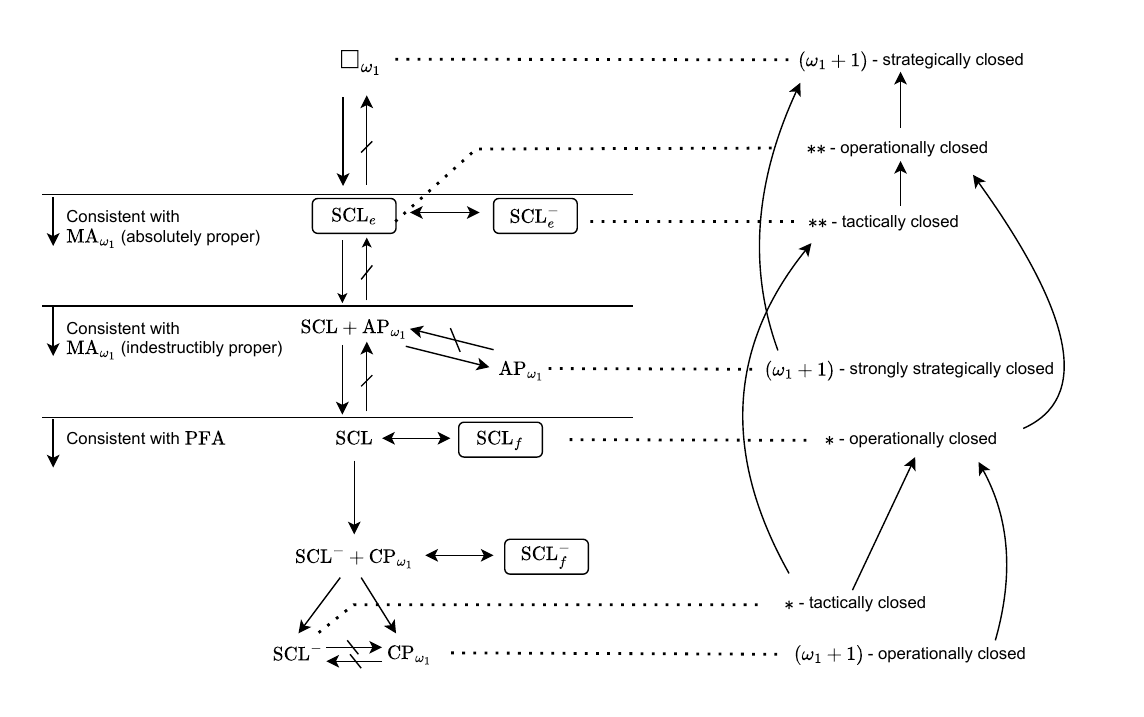}
\end{figure}

%%%%%%%%%%%%%%%%%%%%%%%%%%%%%%%%%%%%%%%%%%%%
\section{Questions}\label{sec:question}
%%%%%%%%%%%%%%%%%%%%%%%%%%%%%%%%%%%%%%%%%%%%
Here we give a list of open questions in this area.

\medskip
Although the statement of $\mathrm{SCL}$ almost looks like a simple combination of $\mathrm{CL}_{\omega_1}$ and $\mathrm{SCL}^-$, we do not know if they are equivalent.
\begin{qtn}\normalfont
Does $\mathrm{CL}_{\omega_1}+\mathrm{SCL}^-$ (that is equivalent to $\mathrm{SCL^-_f}$) imply $\mathrm{SCL}$ ?
\end{qtn}
Theorem \ref{thm:minuse} seems to indicate that the two properties of posets, $**$-tactically closedness and $**$-operationally closedness are close. However, we do not even know whether the conclusion of Theorem \ref{thm:maap} holds for $**$-operationally closed forcing.
\begin{qtn}\normalfont
Is every $**$-operationally closed poset $**$-tactically closed? Or at least does $\mathrm{MA}_{\omega_1}(\text{absolutely proper})$ hold after any $**$-operationally closed forcing, if $\mathrm{PFA}$ holds in the ground model?
\end{qtn}
The fact that $\mathrm{SCL^-_e}$ imply $\mathrm{AP}_{\omega_1}$ seems to indicate that $(\omega_1+1)$-strongly strategic closedness is `stronger' than $**$-tactical closedness in some sense. But we do not know if there is a direct implication.
\begin{qtn}\normalfont
Is every $(\omega_1+1)$-strongly strategically closed poset $**$-tactically closed?
\end{qtn}

%%%%%%%%%%%%%%%%%%%%%%%%%%%%%%%%%%%%%%%%%%%%

%\section*{Acknowledgments}

\bibliographystyle{plain}
\bibliography{locco}
\end{document}